\documentclass[review]{elsarticle}

\usepackage{amssymb,amsmath,comment}
\usepackage{amsthm}
\usepackage[usenames,dvipsnames]{xcolor}
\usepackage{subfigure}
\usepackage{caption}
\usepackage{algorithm, algorithmicx, algcompatible}
\usepackage{algpseudocode}
\usepackage{graphicx}
\usepackage{tikz}
\usepackage[export]{adjustbox}
\usepackage{epstopdf}
\usepackage{hyperref}

\definecolor{cadmiumgreen}{rgb}{0.0, 0.42, 0.24}

\usepackage{lineno}
\modulolinenumbers[5]
\newdefinition{rmk}{Remark }
\newdefinition{prop}{Proposition }

\makeatletter
\def\@author#1{\g@addto@macro\elsauthors{\normalsize%
    \def\baselinestretch{1}%
    \upshape\authorsep#1\unskip\textsuperscript{%
      \ifx\@fnmark\@empty\else\unskip\sep\@fnmark\let\sep=,\fi
      \ifx\@corref\@empty\else\unskip\sep\@corref\let\sep=,\fi
      }%
    \def\authorsep{\unskip,\space}%
    \global\let\@fnmark\@empty
    \global\let\@corref\@empty  %% Added
    \global\let\sep\@empty}%
    \@eadauthor={#1}
}
\makeatother

\begin{document}
\begin{frontmatter}

\title{Fast Converging Parallel Offline-Online Iterative Multiscale Mixed Methods}

\author[]{Dilong Zhou $^1$}
\cortext[cor1]{Corresponding author}
\ead{dxz200000@utdallas.edu}
\author[]{Rafael T. Guiraldello $^2$}
\ead{rafaeltrevisanuto@gmail.com}
\author[]{Felipe Pereira $^1$ \corref{cor1}}
\ead{luisfelipe.pereira@utdallas.edu}

\address{$^1$ Department of Mathematical Sciences, The University of Texas at Dallas, \\
800 W. Campbell Road, Richardson, Texas 75080-3021, USA \\
$^2$ Piri Technologies, LLC \\
 1000 E. University Ave., Dept. 4311, Laramie, WY 82071-2000, USA}


\begin{abstract}
 In this work, we build upon the recently introduced Multiscale Robin Coupled Method with Oversampling and Smoothing (MRCM-OS) to develop two highly efficient iterative multiscale methods. The MRCM-OS methodology demonstrated the ability to achieve flux error magnitudes on the order of $10^{-4}$ in a challenging industry benchmark, namely the SPE10 permeability field. The two newly proposed iterative procedures, through the construction of online informed spaces, significantly enhance the solution accuracy, reaching flux error magnitudes of order $10^{-10}$ for a reduced number of steps.
 
 The proposed methods are based on the construction of online informed spaces, which are iteratively refined to improve solution accuracy. Following an initial offline stage, where known boundary conditions are applied to construct multiscale basis functions, the informed spaces are updated through iterative procedures that utilize boundary conditions defined by the most recently computed solution variables. Two distinct approaches are introduced, each leveraging this framework to deliver efficient and accurate iterative solutions.
 
 A series of numerical simulations, conducted on the SPE10 benchmark, demonstrates the very rapid convergence of the iterative solutions. These results highlight the computational efficiency and competitiveness of the two proposed methods, which are thoroughly compared to each other and to an existing multiscale iterative method from the literature.

\end{abstract}

\end{frontmatter}
%%\linenumbers

\section{Introduction}\label{intro}
Two and three-phase flow problems of interest to the oil industry involve high-resolution permeability data, billions of cells, and long-time steps. Within each time step, a convection-diffusion equation must be solved alongside a second-order elliptic equation. This paper focuses on solving these elliptic equations numerically.

Over the past few decades, multiscale methods have been a key research focus, particularly in addressing complex subsurface flow problems. These approaches generally fall into three main categories: Multiscale Finite Volume, Multiscale Finite Element, and Multiscale Mixed Finite Element methods.

The Multiscale Finite Volume method was developed in 1997 \cite{HouMulti, HouMultivol}, followed by the Multiscale Finite Element Method \cite{jenny2003multi}. This study focuses on Multiscale Mixed Finite Element methods, which are directly related to the work presented here.
The first of these method was introduced in 2003 \cite{chen_hou}, followed by the Multiscale Mortar Mixed Finite Element Method (MMMFEM) in 2007 \cite{arbogast} and the Multiscale Hybrid-Mixed Finite Element Method (MHM) in 2013 \cite{araya2013multiscale}. Our work primarily focuses on Multiscale Mixed Methods with Improved Accuracy: The Role of Oversampling and Smoothing (MRCM-OS) \cite{MRCM-OS}, developed from the Multiscale Robin Coupled Method (MRCM) \cite{guiraldello2018multiscale}, which incorporates two novel strategies: oversampling and smoothing. The MRCM in turn, is a generalization of the Multiscale Mixed Method (MuMM) introduced in 2014 \cite{pereira, hani_pereira_1, 2020MUMM}.

MRCM’s inherent scalability on multi-core supercomputers has made it a valuable tool for solving subsurface flow problems \cite{Recursive2023, HPC2022}. With the help of a post-processing method that ensures continuous normal flux components \cite{velocityMM2020}, MRCM can approximate multiphase flows in porous media with higher accuracy. Moreover, a key parameter, denoted as $\alpha$, controls the relative importance of pressure and the normal flux component in the Robin boundary condition along the boundary. Adjusting this parameter to 0 or $\infty$ allows MRCM to converge the outcomes of the Multiscale Mortar Mixed Finite Element Method (MMMFEM) and the Multiscale Hybrid-Mixed Finite Element Method (MHM), respectively.

Due to MRCM's remarkable scalability on parallel multicore supercomputers, substantial effort has been devoted to enhancing and applying this method to new challenges. See, for instance \cite{interface2021, MMtwo2020, MRCM2022,MPM1,MPM2}.

The MRCM-OS method reduces resonance errors along the interfaces of non-overlapping subdomains compared to MRCM by utilizing oversampling and smoothing procedures. Oversampling allows each Multiscale Basis Function (MBF) to be computed in oversampling regions while maintaining compatibility conditions on the non-overlapping partitions. Smoothing procedures, a multiplicative Schwarz-type method inspired by overlapping Schwarz domain decomposition techniques \cite{DDintro, DomainD}, help correct the small-scale errors inherent in the multiscale solution, further enhancing accuracy. With these two new strategies, MRCM-OS achieves a two-order magnitude improvement compared to MRCM.

However, relying solely on MRCM-OS with linear polynomial spaces, even with the incorporation of oversampling and smoothing, results in an error of approximately $10^{-4}$. Expanding these spaces to achieve smaller errors is not only computationally expensive but also technically challenging, especially in 3D scenarios. To address this, we adopt iterative procedures. The development of iterative methods based on domain decomposition is an active area of research, attracting significant interest from various groups employing diverse approaches. In this context, we introduce a local-global strategy along with an offline-online staged method that utilizes an iterative approach to reduce the error effectively.

The local-global method is primarily used in domain decomposition methods. In each iteration, it processes each subdomain (local) based on the results from the previous iteration or initial conditions for the first iteration. The outcomes are then combined to address the entire domain (global), and the process continues until a specified criterion is met. This method has been utilized in Multiscale methods for some time. In the Multiscale Finite Volume Method \cite{2003FV, 2006FV, 2007FV, 2008FV, 2011FV}, each iteration updates the upscaled transmissibilities $T^*$ or the upscaled permeability $k^*$ based on the previous results. For the Generalized Multiscale Finite Element Method \cite{2005GFE, 2008GFE, 2012GFE, 2018GFE, 2022GFE, 2023GFE, 2024GFE}, the enrichment function is updated based on previous iterations, and the enriched global problem is solved for the current iteration. In the Multiscale Finite Element Method \cite{2011FE, 2014FE, 2017FE}, each iteration updates the Multiscale Basis Function based on the previous iteration and uses Galerkin projection to obtain an approximate solution for the next iteration. The multiscale mortar mixed finite element method taking advantage of a local-global strategy \cite{2019MMMFEM} is used to calculate two-phase flow problems within each time step. 

In the offline-online staged method, part of the calculations can be done with known boundary conditions (offline stages), and some calculations require the result of the offline stage or previous online computations to be performed (Online stage). Typical examples include additive and multiplicative Schwarz methods \cite{2021additive, 2024additive, 2024additive2}, as well as other Schwarz-type methods, such as the smoothing procedure used in \cite{MRCM-OS}, which will also be used here. In the offline stage, an algorithm is used to handle the initial conditions, as demonstrated in \cite{interface2019}, where an interface space is constructed during the offline stage using the Multiscale Robin Coupled Method. During the online stage, which typically involves iteration, the results from previous iterations or the offline stage are used to generate new outcomes. Different types of offline-online staged methods can be explained as follows. The first approach is similar to the local-global method, where both the online and offline stages use similar algorithms, leading to comparable computational costs. This type can be found in \cite{2007same, 2008same, 2013same, 2013sameo, 2013same2, 2014same, 2014more, 2015same, 2017same, 2018same, 2018same2, 2020same, 2022same, 2022same2, 2023same}. The second type involves using a simpler algorithm in the online stage, with lower computational costs compared to the offline stage. This type can be found in \cite{2014same, 2015less}. The third type requires accumulating certain variables during each online stage, leading to increased complexity with more iteration steps. This type can be found in \cite{2010more, 2010more2, 2014more, 2016more}. We introduce two distinct iterative approaches. The first, referred to as the Reduced Method, integrates offline-online stages of types 1 and 2. The second approach, called the Extended Method, encompasses types 2 and 3.

Our current study focuses on enhancing the Multiscale Robin Coupled Method with oversampling and smoothing (MRCM-OS) \cite{MRCM-OS} using the offline-online staged strategy, aiming at minimizing error while maintaining computational efficiency through a tipically minimal number of iteration steps. 
Our main contribution in developing the Reduced and Extended Methods consists of constructing informed spaces along with corresponding offline-online staged methods. In the online stage, the informed spaces generate new Robin boundary conditions for each subdomain based on the pressure and flux results from previous iterations. This approach effectively reduces large-scale errors in the multiscale solution, while MRCM-OS addresses small-scale errors, further enhancing the accuracy of our numerical method.

Through comprehensive studies, we determine the optimal alpha value that converges in the fewest steps and analyze the impact of oversampling and smoothing procedures on the iteration count with this alpha. Our results indicate that, with the optimal alpha, our offline-online method achieves flux accuracy of $10^{-10}$ within 10 iterations without requiring extensive oversampling or numerous smoothing steps, even for challenging permeability fields.

This paper is organized into several sections as follows. The first section focuses on the formulation of MRCM with oversampling. Next, our offline-online staged methods are described, with the Reduced and Extended Methods each detailed in their own subsection. The subsequent section presents our numerical studies, subdivided into subsections to explore various aspects of our experiments. We consider problems involving two permeability fields derived from the SPE 10 project \cite{SPE10}, one with a challenging channel and the other with the highest contrast ratio. At the end of this section, we compare our method with an existing method \cite{2008FV} on the same problem. Through these numerical studies, we demonstrate the effectiveness and improved performance of our MRCM-OS combined with the offline-online staged strategy.

\section{The MRCM with Oversampling (MRCM-O)}

We consider single-phase flow in porous media. The governing equations for pressure $p$ and Darcy velocity ${\bf u}$ are given by
\begin{eqnarray}
{\bf u} & = & -\, K\,\nabla p \qquad \mbox{in}~\Omega \label{eq1a}\\
\nabla \cdot {\bf u} & = & f \qquad \mbox{in}~\Omega \label{eq1b}\\
p & = & g \qquad \mbox{on}~\partial \Omega_p \\
{\bf u}\cdot {\bf n} & = & z \qquad \mbox{on}~\partial \Omega_u \label{eq1d}
\end{eqnarray}
where $\Omega\subset \mathbb{R}^d$, $d=2$ or $3$ is the
domain of the problem, $K$ is a symmetric, uniformly positive definite tensor
with components in $L^\infty(\Omega)$,
$f\,\in\,L^2(\Omega)$ is the source term, $g\,\in\,H^{\frac12}(\partial\Omega_p)$
is the pressure condition on the boundary, $z\,\in\,H^{-\frac12}(\partial\Omega_u)$
is the normal velocity condition on the boundary and ${\bf n}$ is the outer normal to $\partial{\Omega}$.

%%%%HERE %%%%%
\subsection{MRCM-O}

\begin{figure}[H]
    \centering
    \includegraphics[width = 0.4\textwidth]{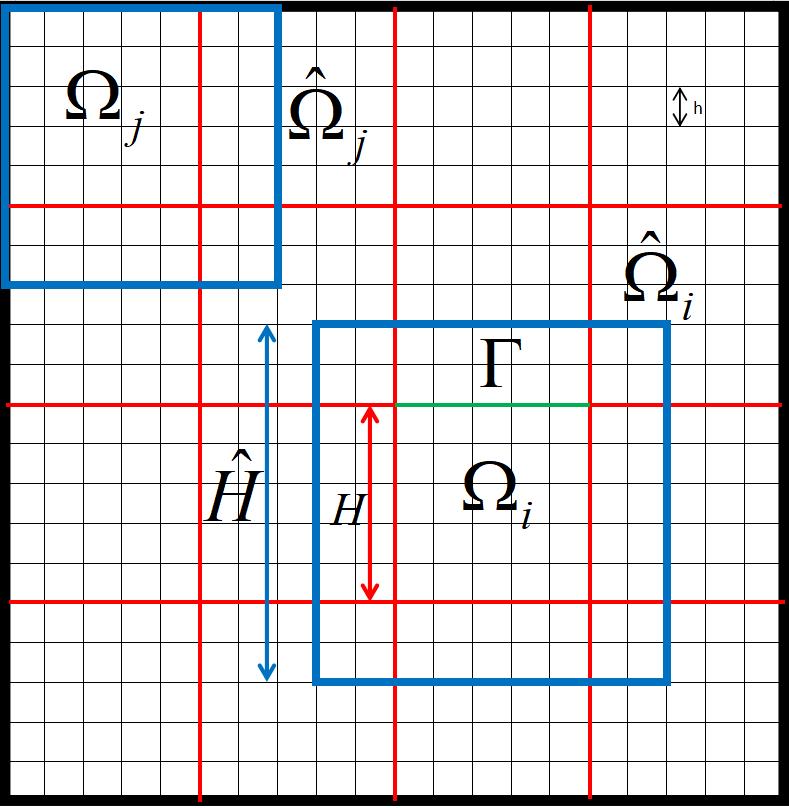}
    \caption{Decompositions of the computational domain $\Omega$: The non-overlapping subdomains are denoted by $\Omega_i$, and the corresponding oversampling regions are written as $\hat{\Omega}_i$. The three length scales that enter in the formulation of the MRCM with oversampling are also illustrated:  $\hat H > H \geq h$.}
    \label{overlapjpg}
\end{figure}

The approach utilized in this study is founded on the MRCM-O \cite{MRCM-OS} which is briefly summarized in this section. Consider $\mathcal{T}_h$ as a subdivision of $\Omega\subset \mathbb{R}^d$ into a Cartesian mesh of $d$-dimensional rectangles. From this mesh, define partitions of $\Omega$ into non-overlapping subdomains $\{\Omega_i\}_{i=1,\ldots,m}$, and let $\Gamma$ be the union of all interfaces $\Gamma_{i,j}=\overline{\Omega}_i\cap\overline{\Omega}_j, \forall i,j=1\dots\,m$, and define $\Gamma_i = \partial\Omega_i\setminus\partial\Omega$. For each subdomain $\Omega_i$, define $\hat{\Omega}_i$ as the augmented subdomain comprising $\Omega_i$ along with an adjoining region. Figure \ref{overlapjpg} illustrates these definitions.
 The discrete variational formulation of the Multiscale Robin Coupled Method with Oversampling (MRCM-O) reads as: 
Find $({\bf u}_h^i,p_h^i,\lambda_h^i)\,\in\,{\bf V}_{hz}^i\times Q_h^i
\times \Lambda^i_{H}$, for
$i=1,\ldots,m$, such that
\begin{eqnarray}
  (K^{-1}{\bf u}_h^i,{\bf v})_{\Omega_i}-(p_h^i,\nabla\cdot {\bf v})_{\Omega_i} 
  +(\beta_i\,{\bf u}_h^i\cdot\check{\bf n}^i,{\bf v}\cdot\check{\bf n}^i)_{\Gamma_i}
  +(\lambda_h^i,{\bf v}\cdot\check{\bf n}^i)_{\Gamma_i} \nonumber \\ = -(g,{\bf v}\cdot\check{\bf n}^i)_{\partial\Omega_i\cap\partial\Omega_p},  \label{eq11d}\\
  (q,\nabla\cdot {\bf u}_h^i)_{\Omega_i}  =   (f,q)_{\Omega_i}, \label{eq12d}\\
  \sum_{i=1}^m ({\bf u}_h^i\cdot \check{\bf n}^i,M)_{\Gamma_i} = 0,  \label{eq13d}\\
  \sum_{i=1}^m (\beta_i\,{\bf u}_h^i\cdot \check{\bf n}^i+\lambda_h^i,V \,\check{\bf n}^i\cdot\check{\bf n})_{\Gamma_i} = 0,  \label{eq14d}
\end{eqnarray}
hold for all $({\bf v},q)\,\in\,{\bf V}_{h0}^i$ and for all $(M,V)\,\in\, M_H \times V_H \subset F_h(\mathcal{E}_h)\times\,F_h(\mathcal{E}_h)$, where
\begin{eqnarray*}
  {\bf V}_{hy}^i&=&\{{\bf v}\,\in\,{\bf V}_{h}^i~,~
  {\bf v}\cdot \check{\bf n}=y~\mbox{on}\,\partial\Omega_i\cap\partial\Omega_u
  \}~, \label{eq:Vhy}
\end{eqnarray*}
and ${\bf V}_{h}^i\times Q_h^i$ are the lowest-order Raviart-Thomas \cite{RaviartThomas::1977}
spaces for velocity and pressure defined for $\Omega_i$, and
\begin{equation*}
  F_h(S_h) = \{ f:{S}_h\to \mathbb{R}~|~f|_e\,\in\,\mathbb{P}_0~,
  ~\forall\,e\,\in\,{S}_h \}
\end{equation*}
with $\mathcal{E}_h$ defined as the set of all faces/edges of $\mathcal{T}_h$ contained 
in $\Gamma$.

The Lagrange multiplier spaces $\Lambda^i_H = \mbox{span}\left\{\phi_i^1,\phi_i^2,..,\phi_i^N\right\} \subset F_h(\mathcal{E}_h\cap\Gamma_i)$ with $\phi_i^k$ given by
\begin{equation*}
  \phi_i^k = -\beta_i\,{\bf u}_h^{k}\cdot {\bf \check{n}}^i|_{\Gamma_i} + \pi^k|_{\Gamma_i}, \ \forall\,k=1,\dots,N,  
  \label{informed_space}  
\end{equation*}
where ${\bf u}_h^{k}\cdot {\bf \check{n}}^i$ and $\pi^k$ is the normal component of velocity
and interface pressure, respectively, retrieved on interface $\Gamma_i$ as solutions of the discrete
problem given by
\begin{equation}
  \begin{array}{rclll}
    {\bf u}_h^{k} &=& -K\nabla_h\,p_h^{k} &&\mbox{in} \ \hat\Omega_i \\
    \nabla_h\cdot{\bf u}_h^{k} &=& 0 &&\mbox{in}  \ \hat\Omega_i \\
    p_h^{k} & = & 0 \qquad &&\mbox{on}~\partial\hat\Omega_i\cap\partial\Omega_p \\
    {\bf u}_h^{k}\cdot {\bf n}^i & = & 0 \qquad &&\mbox{on}~\partial\hat\Omega_i\cap\partial\Omega_u \\
    -\beta_i\,{\bf u}_h^{k}\cdot {\bf n}^i + p_h^{k} &=& \lambda^k &&\mbox{on}~\partial\hat\Omega_i\setminus\partial\Omega
  \end{array}, \label{eq:oversampling_problem}
\end{equation}

To efficiently solve the system (\ref{eq11d})--(\ref{eq14d}) while leveraging the computational benefits of multi-core architectures, we adopt the concept of multiscale basis functions (\emph{MBFs}) (see \cite{yotov2009}, \cite{guiraldello2018multiscale}, \cite{MRCM-OS}). Upon solving equations (\ref{eq11d})--(\ref{eq14d}), a smoothing scheme is applied for a specified number of iterations (details of the smoothing procedure can be found in \cite{MRCM-OS}). This comprehensive approach is referred to as \emph{MRCM-OS}, and it is outlined in the algorithm \ref{alg:mrcm_os}.

\begin{algorithm}
  \caption{Algorithm describing the computation of the MRCM-OS solution.}
  \label{alg:mrcm_os}
  \begin{algorithmic}[1]
    \State \textbf{Input: }$\text{Domain\_decomposition}, \text{oversampling\_size}$
    \State \textbf{Step 1} Construct the Lagrange multiplier spaces $\Lambda^i_H$ by solving \eqref{eq:oversampling_problem}.
    \State \textbf{Step 2} Solve the system (\ref{eq11d})--(\ref{eq14d}) trough MBFs.
    \State \textbf{Step 3} Apply the smoothing scheme iteratively for the specified number of steps.
  \end{algorithmic}
\end{algorithm}

A comprehensive description of the MRCM-OS can be found in \cite{MRCM-OS}.

\section{Iterative MRCM with Oversampling and Smoothing (i-MRCM-OS)}

In this section, we describe two iterative algorithms designed to enhance the MRCM-OS solution by incorporating the construction of new MBFs or online multiscale basis functions based on previous MRCM-OS solution. These algorithms aim to efficiently adapt and refine the multiscale solution during the simulation process, ensuring improved accuracy.

\subsection{Reduced Method}

The first methodology involves the iterative construction of new multiscale basis functions derived from the previous MRCM-OS solutions, which then replace the existing set of basis functions. In each iteration, a new set of basis functions is built to solve the MRCM-OS method while preserving the number of degrees of freedom of the initial problem. The algorithm describing the method is shown in Algorithm \ref{method_1}.

\begin{algorithm}[h]
  \caption{Reduced Method}
  \begin{algorithmic}[1]
    \State \textbf{Initialization}:
    \State Compute MBFs with piecewise constant boundary conditions.
    \State Solve MRCM-OS with piecewise constant interface functions.
    \For{$i=1$ to $M$}
    \State Update the MBFs using the revised Robin boundary condition from the previous solution.
    \State Solve MRCM-OS with piecewise constant interface functions.
    \EndFor
  \end{algorithmic}\label{method_1}
\end{algorithm}

The process begins by constructing multiscale basis functions that use piecewise constant boundary conditions, i.e., the boundary conditions of Eq.~\eqref{eq:oversampling_problem} are chosen as 
piecewise constant. This initial choice works as a starting guess in the iterative scheme.
Following the construction of the basis functions, the subsequent step entails solving the MRCM-OS problem with piecewise constant test functions, where the coarse space $M_H \times V_H$ in formulation \eqref{eq11d}-\eqref{eq14d} is spanned by piecewise constant basis functions. 

The method proceeds into a looping for a number of $M$ steps defined by the user and works as follows:
New $\lambda^k$ functions are constructed using the pressure and flux computed from the previous iteration to define the new Robin boundary conditions in Eq.~\eqref{eq:oversampling_problem} and compute the new set of multiscale basis functions. With the updated basis functions, the MRCM-OS is solved again with piecewise constant functions for the interface space. 

The process of updating the basis functions, solving the MRCM-OS continues for the desired number of steps.

\subsection{Extended Method}
The second methodology involves constructing an augmented set of multiscale basis functions. In this approach, the initial basis functions, computed with piecewise constants, are preserved. Based on the current MRCM-OS solution using these initial basis functions, a set of informed multiscale basis functions is computed. This  new set augments the Lagrange multiplier space.  

The method proceeds into a looping for a number of $M$ steps defined by the user and works as follows: 
In each iterative step, based on the current set of multiscale basis functions, the MRCM-OS problem is solved. A new set of informed multiscale basis functions is constructed and replaces
the previous informed multiscale basis functions. The algorithm describing the method is shown in Algorithm \ref{method_3}.

\begin{algorithm}[htb]
  \caption{Extended Method}
  \begin{algorithmic}[1]
    \State \textbf{Initialization}:
    \State Compute MBFs with piecewise constant boundary conditions.
    \State Solve MRCM-OS with piecewise constant interface functions.
    \State Compute informed-MBFs based on current solution.
    \State Append informed-MBFs to the initial MBFs set. (augmented set)
    \For{$i=1$ to $M$}
    \State Solve MRCM-OS with piecewise linear interface functions.
    \State Compute informed-MBFs based on current iteration.
    \State Update the informed-MBFs in the augmented set.
    \EndFor
  \end{algorithmic}\label{method_3}
\end{algorithm}

Note that as the set of multiscale basis functions is augmented, the number of interface functions must be increased accordingly to maintain consistency in the methodology. This is why in Step 3 of Algorithm \ref{method_3}, piecewise constant interface functions are used, and in Step 7, piecewise linear functions are considered.

We note that in this study, the number of informed MBFs is chosen to match exactly the number of MBFs constructed with piecewise constant functions. This ensures that a linear interface maintains consistency in our methodology. A different number of informed MBFs can be selected; however, to maintain consistency, \textbf{Proposition 1.} in \cite{MRCM-OS} should be followed.

It should be noted that RM and EM compute the same number of Multiscale Basis Functions per iteration. However, EM is more computationally demanding because it involves solving an interface linear system that is twice as large as the one required in RM.

\section{Numerical Studies}
The simulations were conducted on an HPC cluster, focusing on the numerical solutions of a model problem involving a highly heterogeneous permeability field from the SPE10 project \cite{SPE10} (http://www.spe.org/web/csp).

\begin{figure}[H]
    \centering
    \includegraphics[width = 0.48\textwidth]{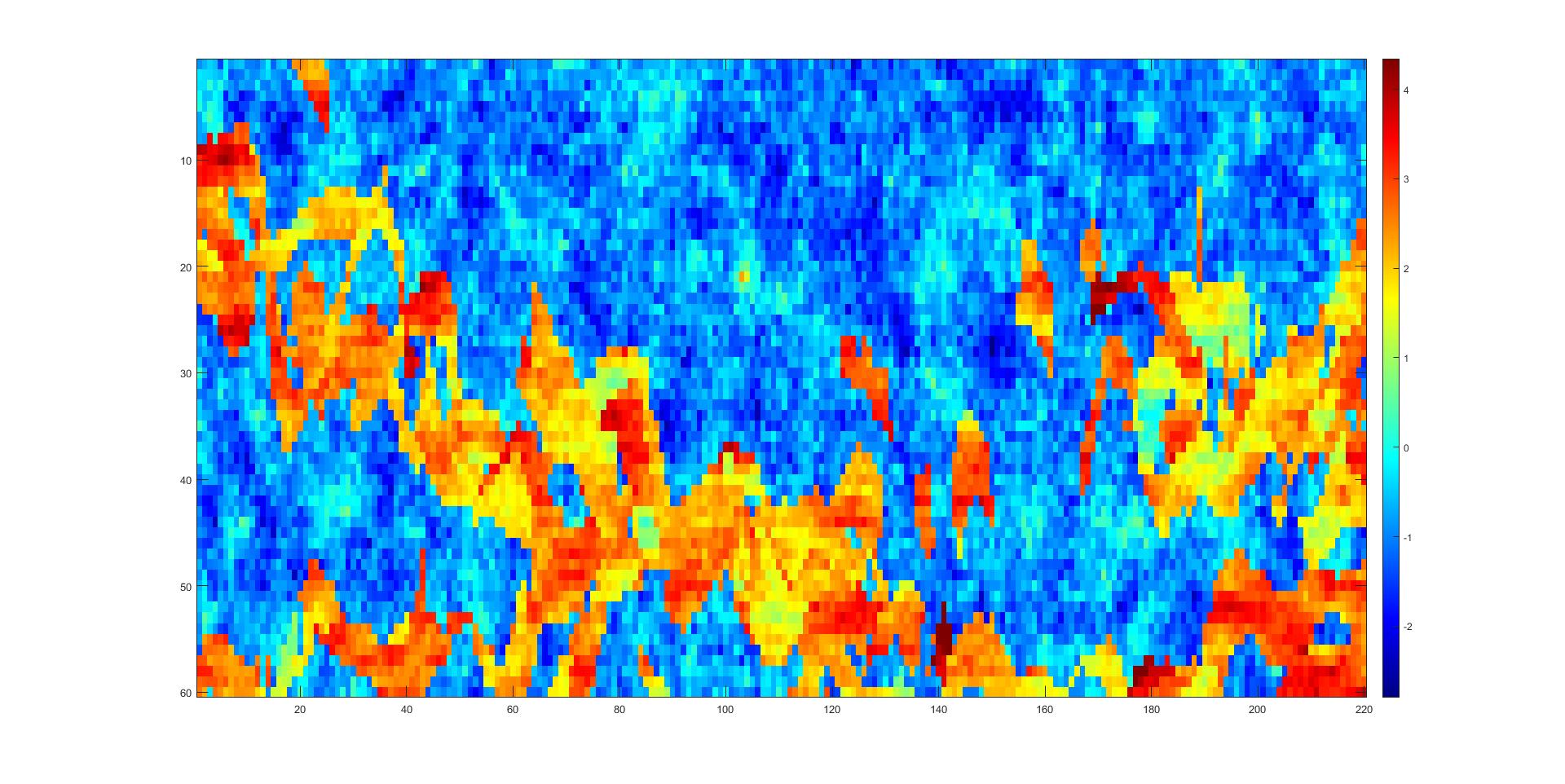}
    \includegraphics[width = 0.48\textwidth]{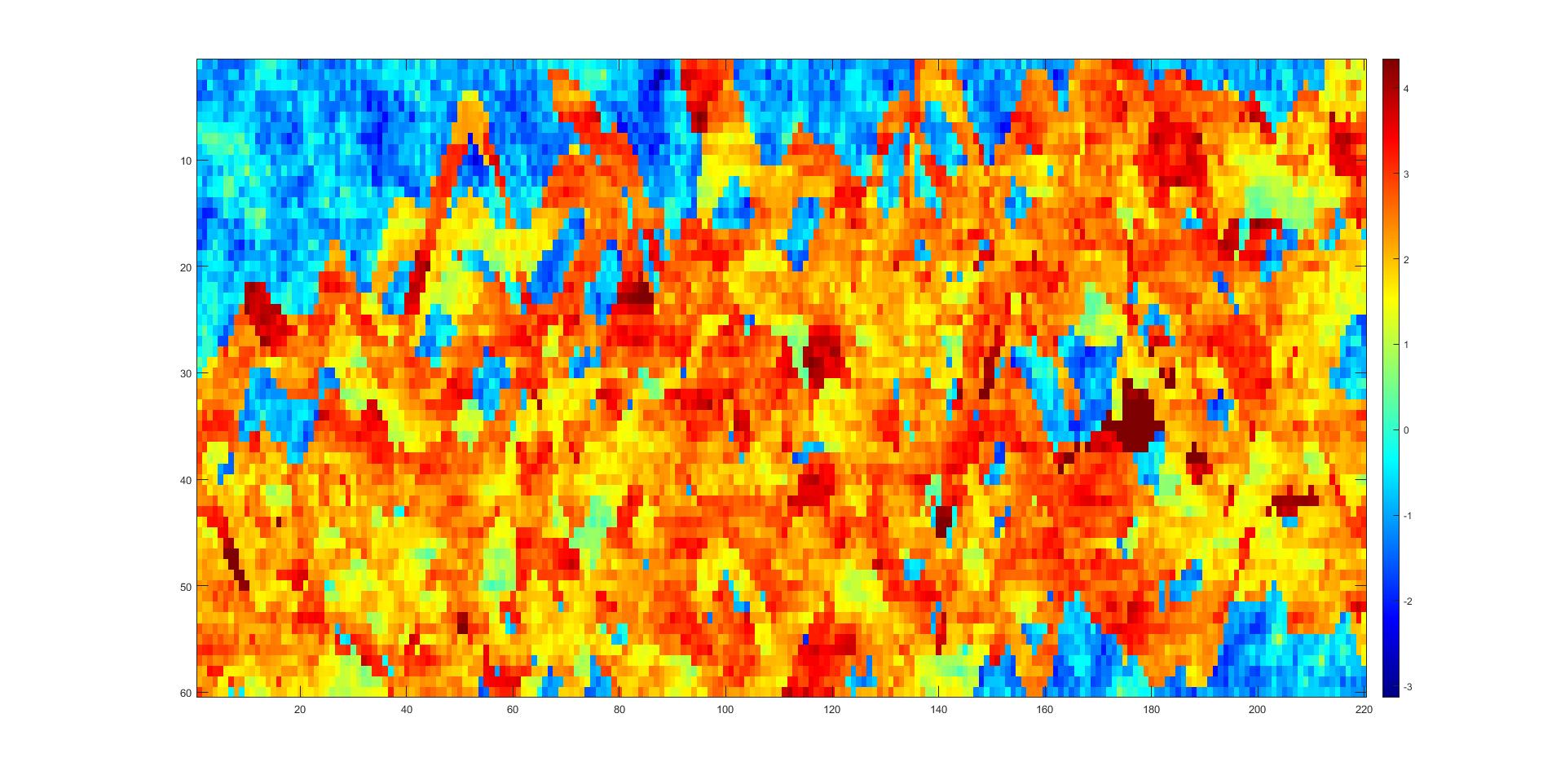}
    \caption{Permeability field: Slice 42 (left) and slice 84 (right) of the SPE 10 project.}
    \label{perm}
\end{figure}

We consider a rectangular domain with Dirichlet boundary conditions, where the pressure is fixed at 1 on the left boundary and 0 on the right boundary. Neumann boundary conditions are imposed with a zero flux at both the top and bottom boundaries, and the source term is set as $q(x)=0$. Additionally, the permeability field $K(x)$ is set to be the 42th and 84th (two-dimensional) layer of the SPE10 project, as depicted in Fig. \ref{perm}. Slice 42 poses a considerable challenge for numerical solvers due to the variability in the permeability field and the presence of a highly permeability channel, which significantly affects the solution behavior and increases the complexity of accurately capturing the flow dynamics. Slice 84, despite lacking an obvious channel, exhibits the highest permeability contrast, with the largest ratio between maximum and minimum permeability values among all slices in the SPE10 project. For all reported results in heterogeneous problems, the domain is decomposed into $11 \times 3$ subdomains with  $H=1/3$, each containing a local $20 \times 20$ grid.

 Table \ref{tab:notation} summarizes the notation to be considered in the numerical experiments. All relative errors in section 4.1 through 4.4 are computed using the  $L^2(\Omega)$ norm, comparing the pressure and flux variables to the fine-grid solution.

\begin{table}[htbp]
  \centering
  \begin{tabular}{|l|l|}
    \hline
    \textbf{Notation} & \textbf{Description} \\
    \hline
    \textit{RM} & Reduced Method \\
    \textit{EM} & Extended Method \\
    \textit{$lh, kS$} & Solution with oversampling size $lh$,\\
    \textit{} & followed by $k$ smoothing steps \\
    \hline
  \end{tabular}
  \caption{Notation used in numerical experiments.}
  \label{tab:notation}
\end{table}

Section 4.1 presents the numerical results from the alpha study for both slices, identifying the optimal alpha values for i-MRCM-O that reduce the relative error to below $10^{-7}$. We continue using this optimal alpha for the subsequent sections. Section 4.2 demonstrates how quickly convergence can be achieved under reasonable settings for oversampling size and smoothing steps. Section 4.3 explores the benefits introduced by oversampling, while Section 4.4 highlights the advantages provided by smoothing. Section 4.5 will present a comparison of our procedures with an existing method \cite{2008FV}
 .
\subsection{Alpha Parameter Studies}
\begin{figure}[H]
    \centering
    \includegraphics[width = 0.48\textwidth]{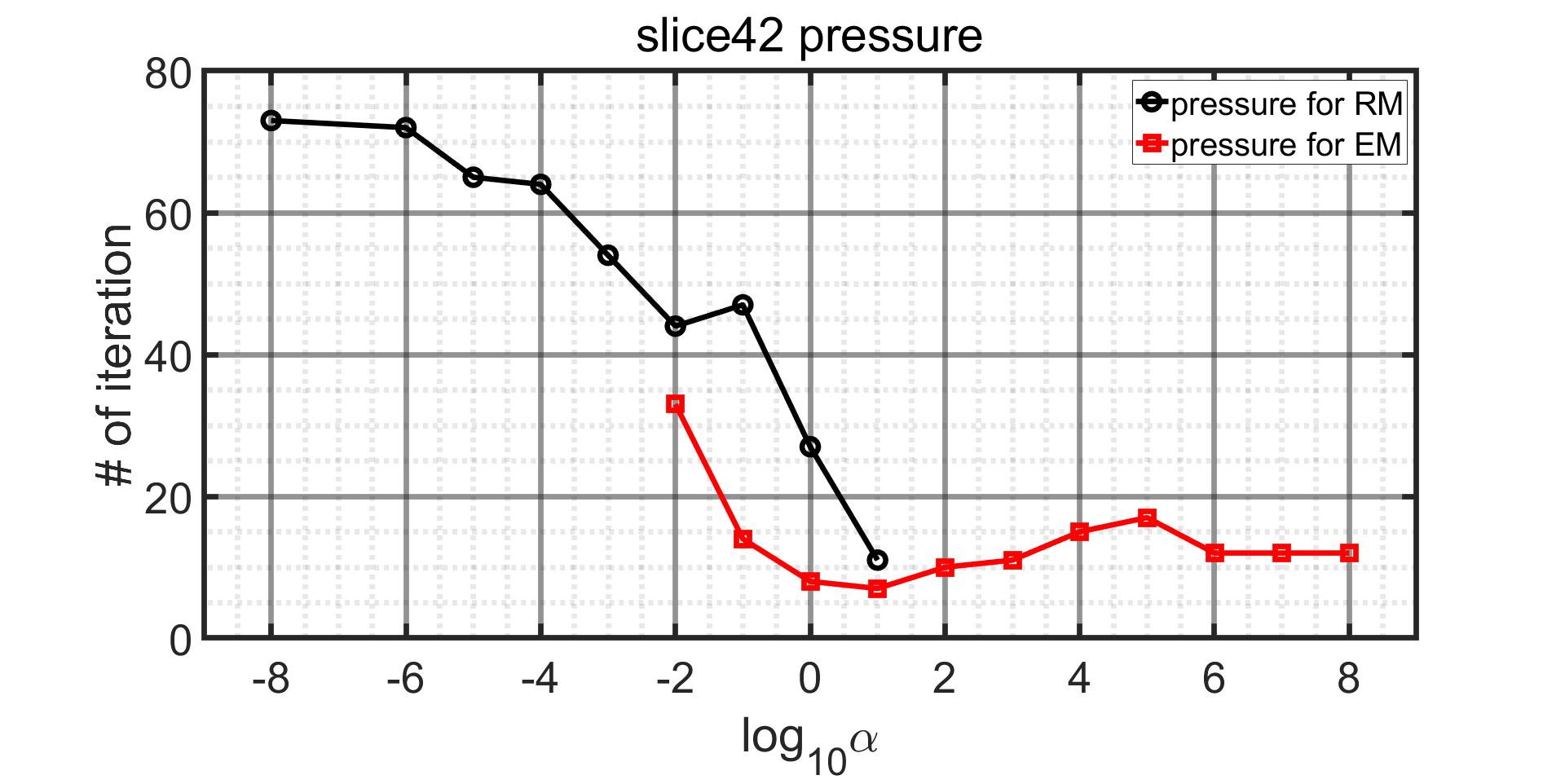}
    \includegraphics[width = 0.48\textwidth]{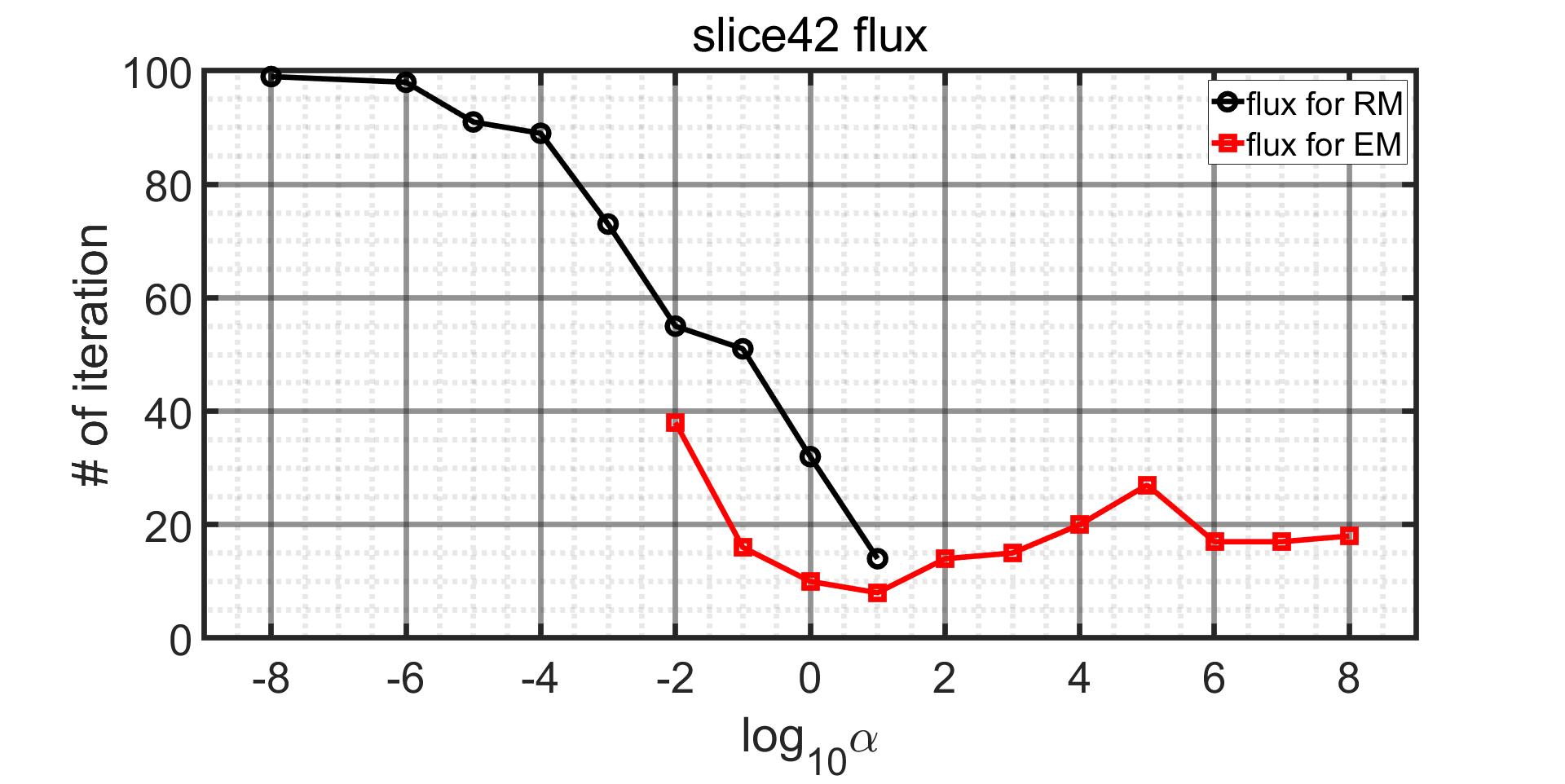}
    \caption{$\alpha$ parameter study for slice 42: iteration times for pressure (left) and flux (right)}
    \label{alpha42}
\end{figure}

\begin{figure}[H]
    \centering
    \includegraphics[width = 0.48\textwidth]{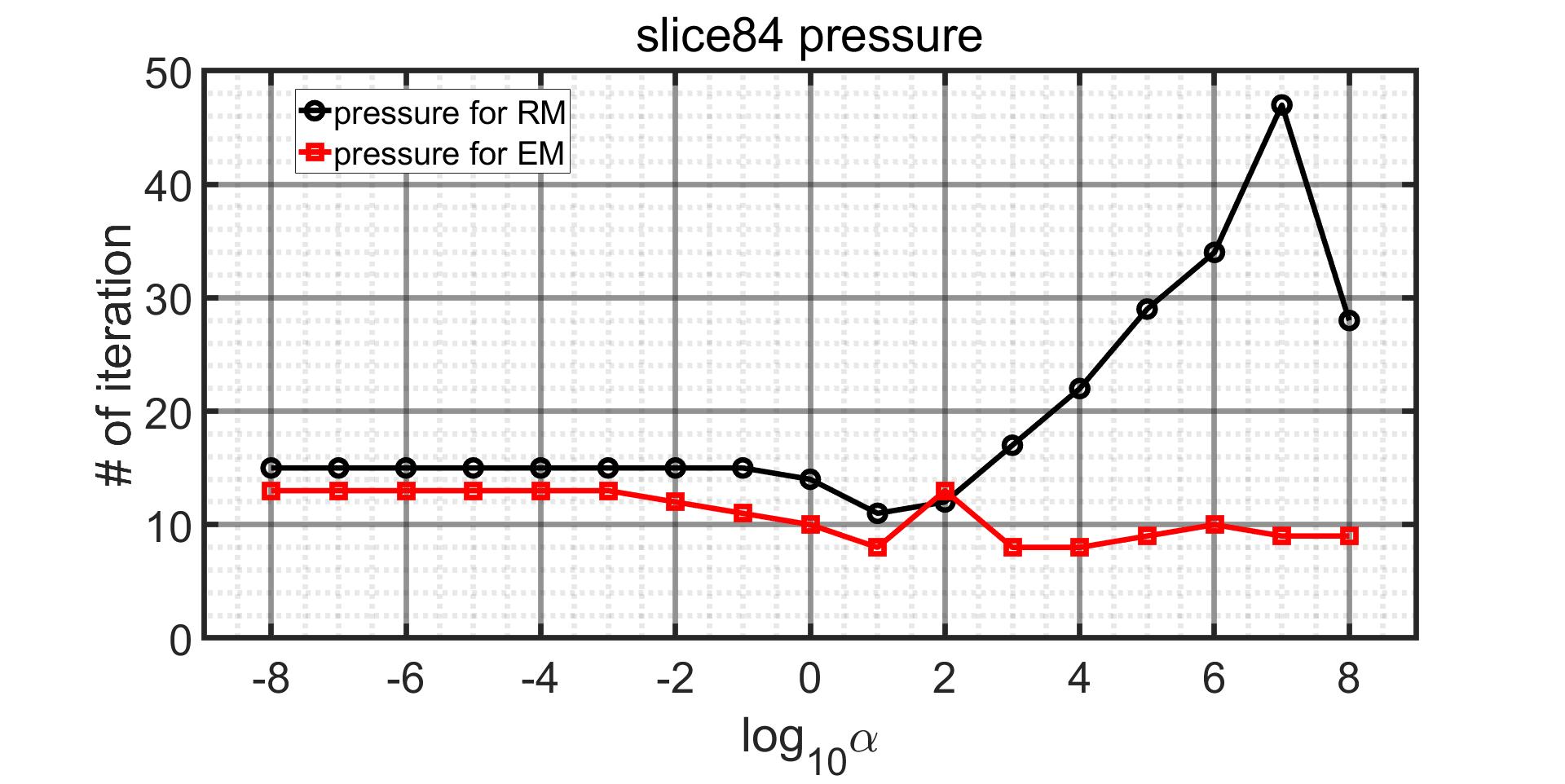}
    \includegraphics[width = 0.48\textwidth]{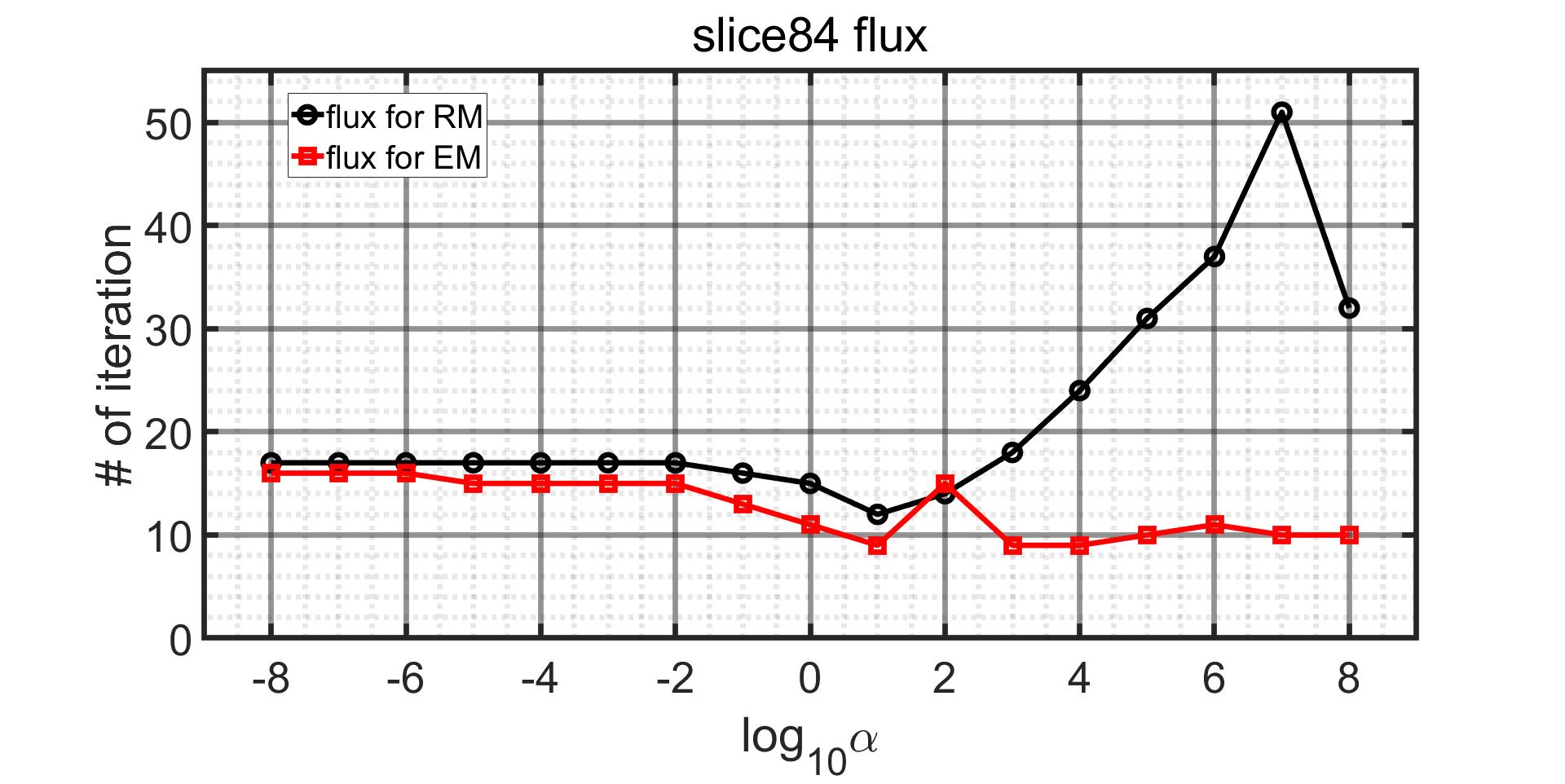}
    \caption{$\alpha$ parameter study for slice 84: iteration times for pressure (left) and flux (right)}
    \label{alpha84}
\end{figure}

Our initial findings are presented in Fig. \ref{alpha42} and \ref{alpha84}. This alpha study is computed for $alpha= 10^{-8}, 10^{-7},...,10^{7},10^{8}$. This study demonstrates the number of iterations required for different alpha values, under the conditions of an oversampling size of 2h and 4 smoothing steps, to reduce the relative error in pressure or flux variables to below $10^{-7}$. A smaller number indicates faster convergence. If a specific alpha lacks a corresponding point in the figure, it means that for that alpha, the relative error did not reach $10^{-7}$ within 100 iterations. The reason could be that it either requires more than 100 iterations or becomes stuck at an intermediate value, such as $10^{-5}$, or oscillates around a value significantly higher than $10^{-7}$.

A key observation from our analysis is that smaller values of $\alpha$ (resembling MMMFEM) and larger $\alpha$ (akin to MHM) influence the Reduced Method and Extended Method differently in these two slices. In slice 84, the Reduced Method performs better with smaller $\alpha$, while the Extended Method shows better results with larger $\alpha$. Notably, across all tested values of $\alpha$, $\alpha = 10$ consistently requires the fewest iteration steps to meet the criteria.

Consequently, we will use only $\alpha = 10$ in the following subsection.

\subsection{Convergence Analysis}
In this section, we highlight results that show the rapid convergence of our method to a very small error.

Since the Extended Method consistently converges faster than the Reduced Method, albeit with higher computational cost, this section focuses solely on the Extended Method. Two scenarios are considered: one with an oversampling size of 2h and 4 smoothing steps (as discussed in Section 4.1), and another with an oversampling size of 4h and 8 smoothing steps.

\begin{figure}[H]
    \centering
    \includegraphics[width = 0.48\textwidth]{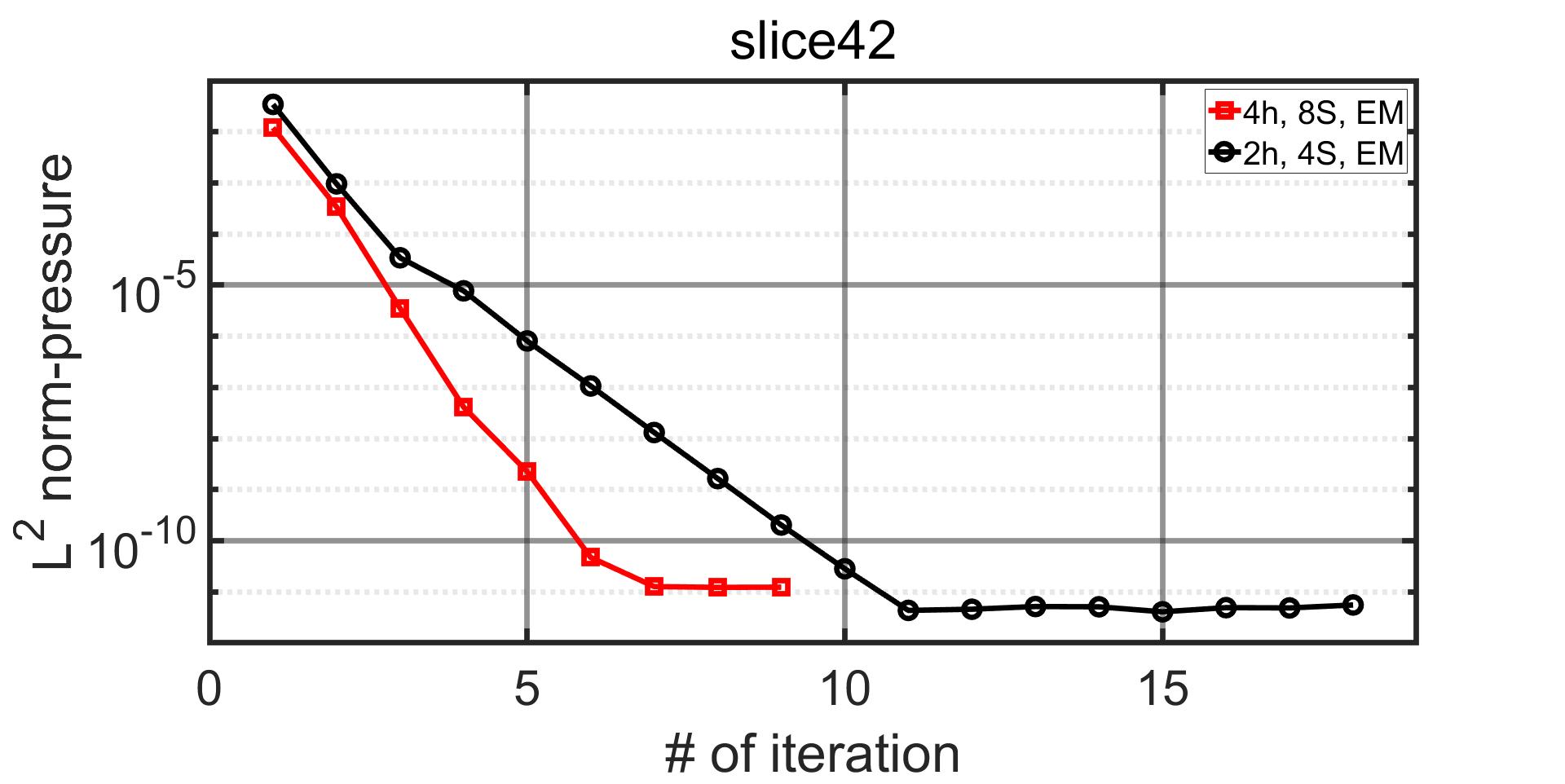}
    \includegraphics[width = 0.48\textwidth]{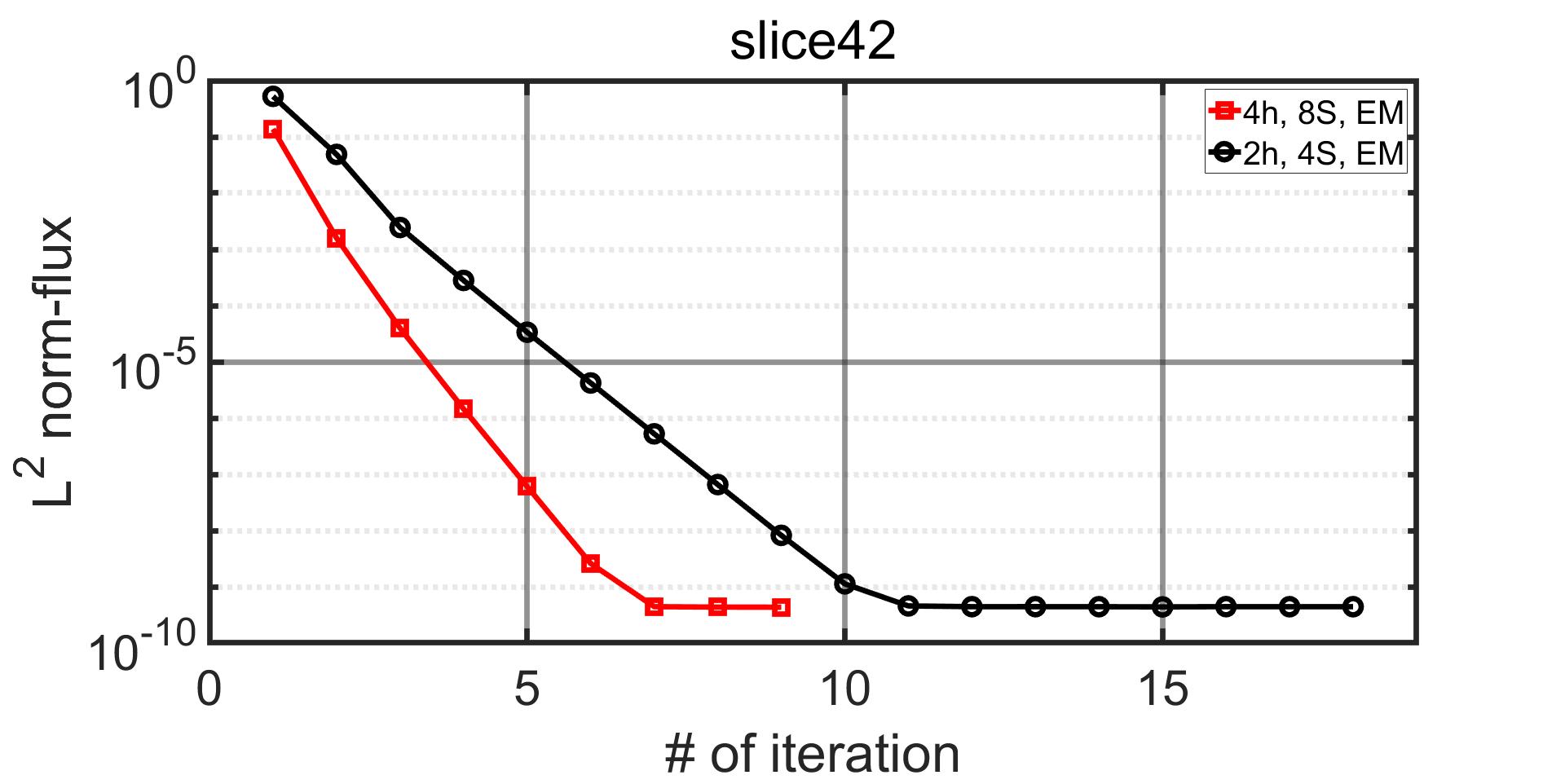}
    \caption{Number of iterations until convergence for slice 42: pressure (left) and flux (right).}
    \label{best42}
\end{figure}

\begin{figure}[H]
    \centering
    \includegraphics[width = 0.48\textwidth]{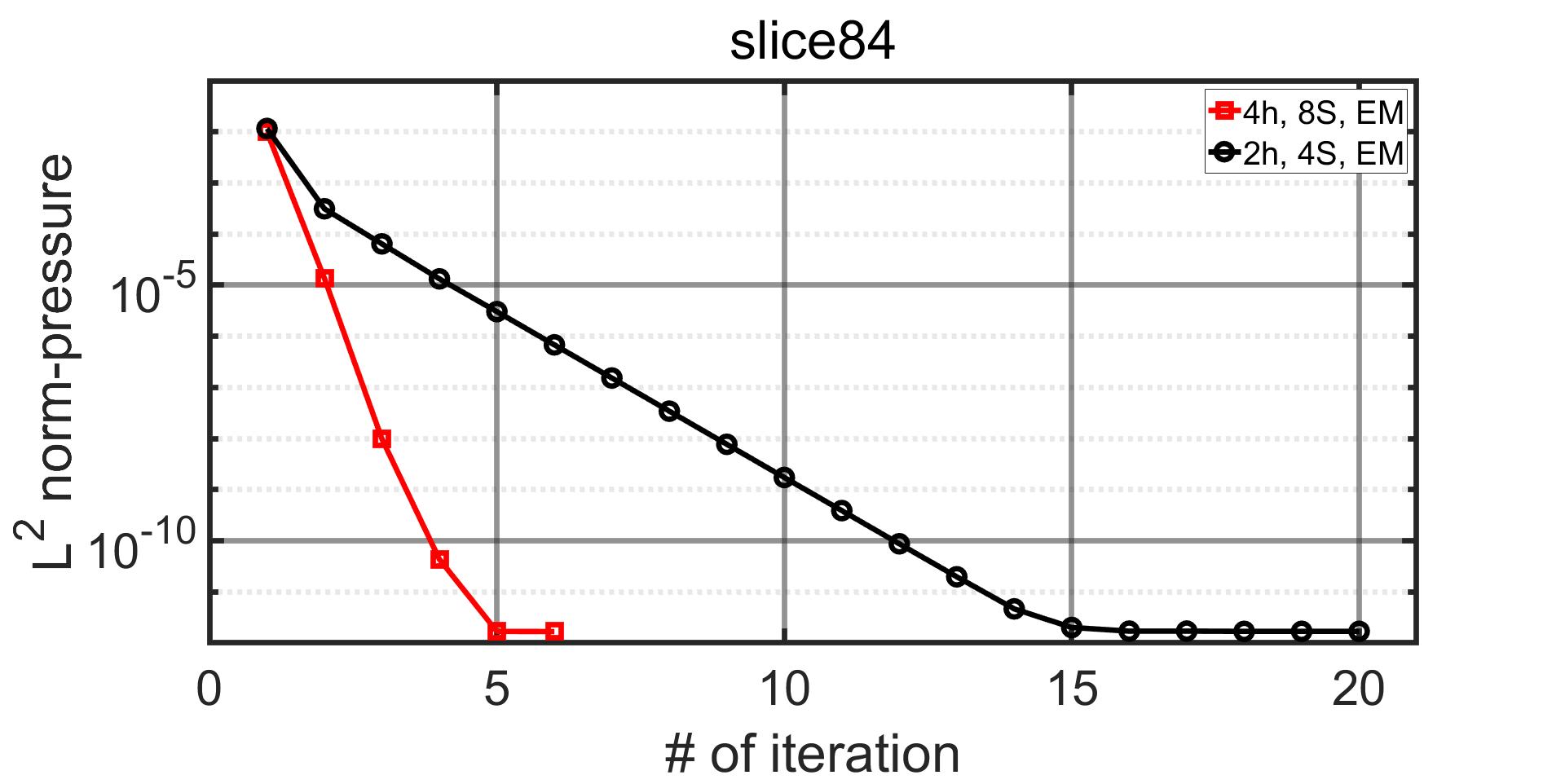}
    \includegraphics[width = 0.48\textwidth]{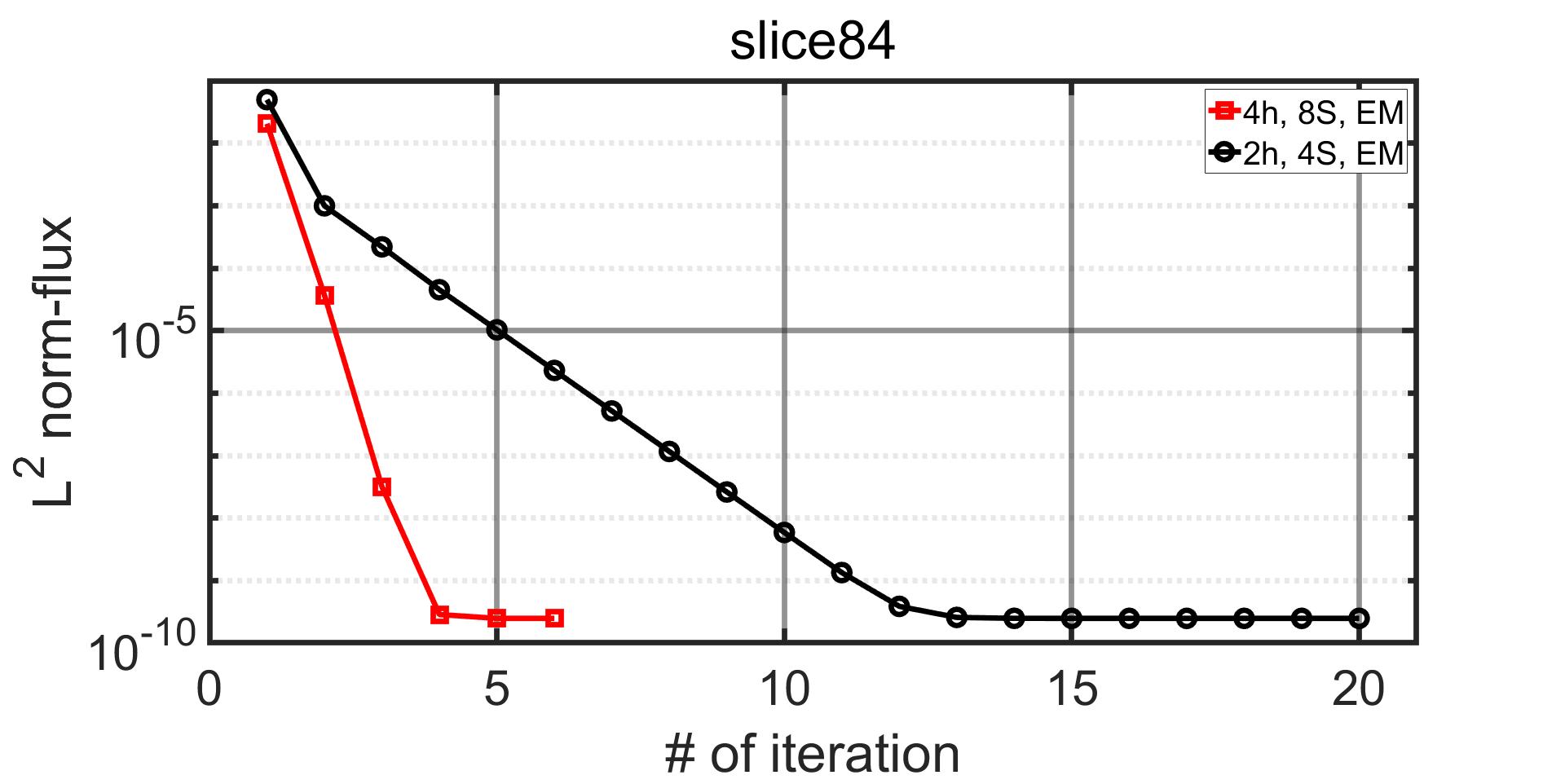}
    \caption{Number of iterations until convergence for slice 84: pressure (left) and flux (right).}
    \label{best84}
\end{figure}

As demonstrated in Figures \ref{best42} and \ref{best84}, under optimal conditions, the method converges rapidly to a very small error within just 7 steps.

In the following subsection, we will discuss how the oversampling size and the number of smoothing steps influence the number of iterations required for convergence.

\subsection{The Role Of Oversampling}
This section presents the results on how oversampling impacts the number of iterations needed for convergence, given fixed smoothing steps.

\begin{figure}[H]
    \centering
    \includegraphics[width = 0.48\textwidth]{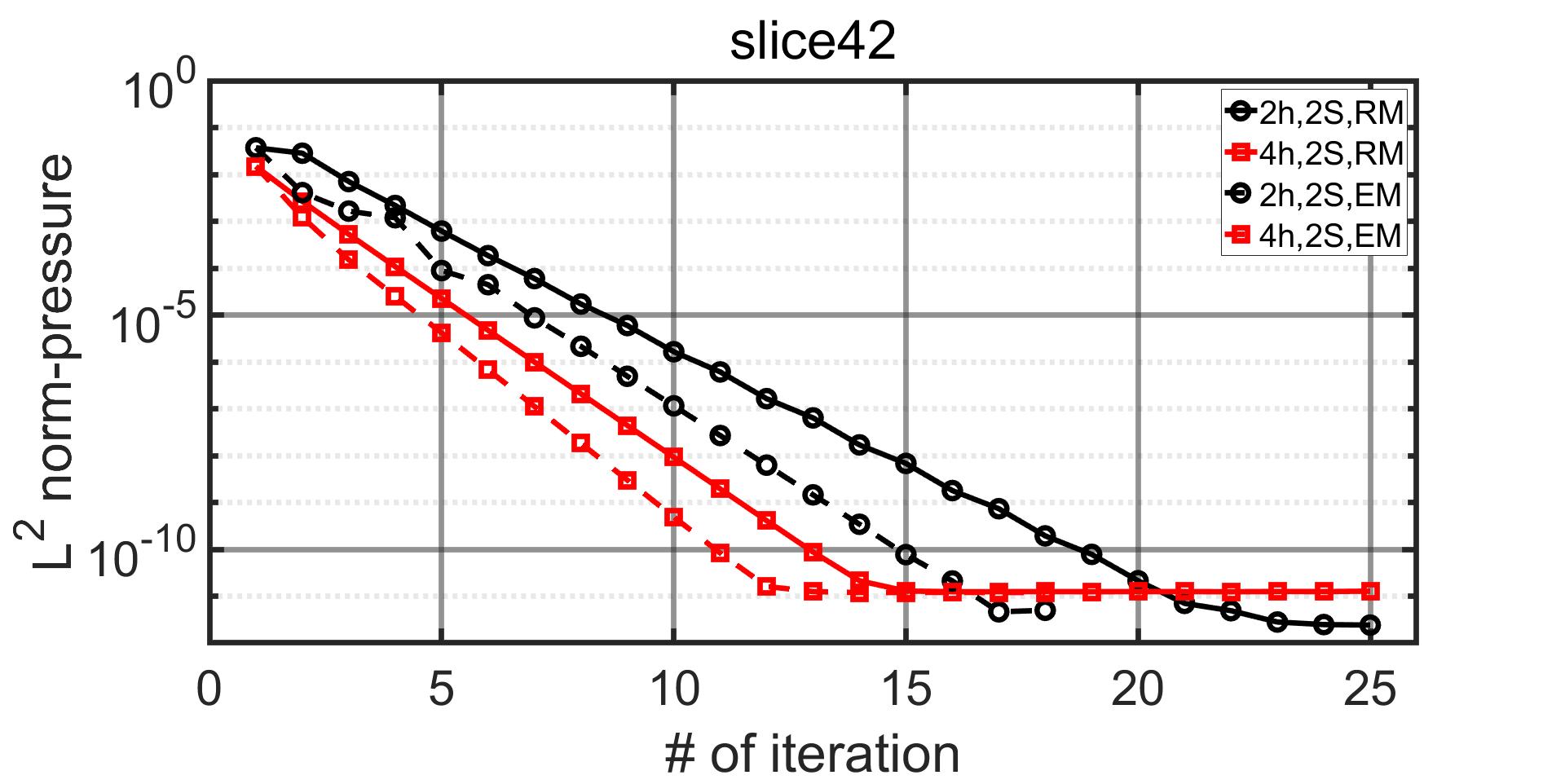}
    \includegraphics[width = 0.48\textwidth]{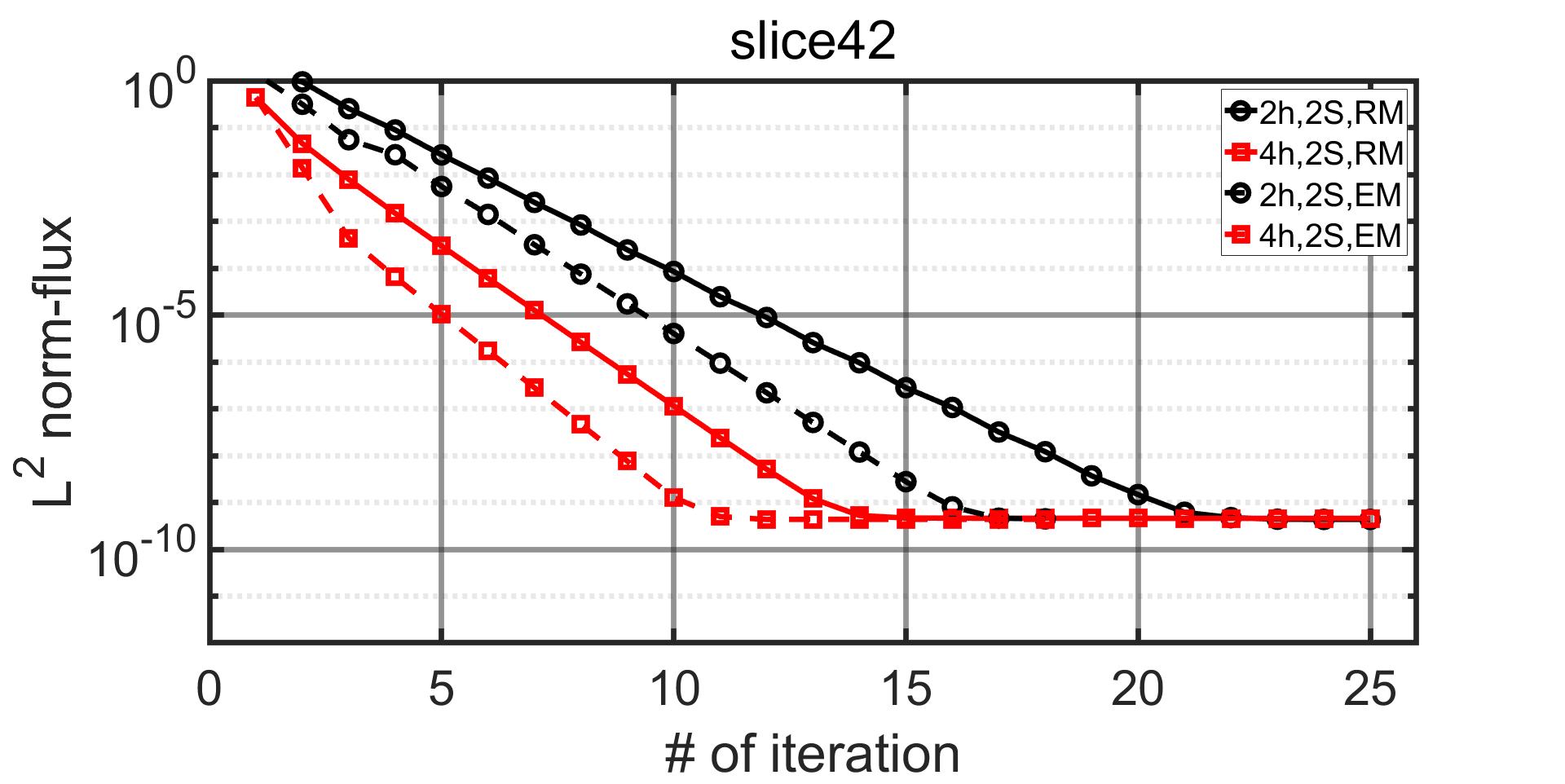}
    \caption{Number of iterations to reach convergence for different oversampling sizes in slice 42: pressure (left) and flux (right).}
    \label{overlap42}
\end{figure}

\begin{figure}[H]
    \centering
    \includegraphics[width = 0.48\textwidth]{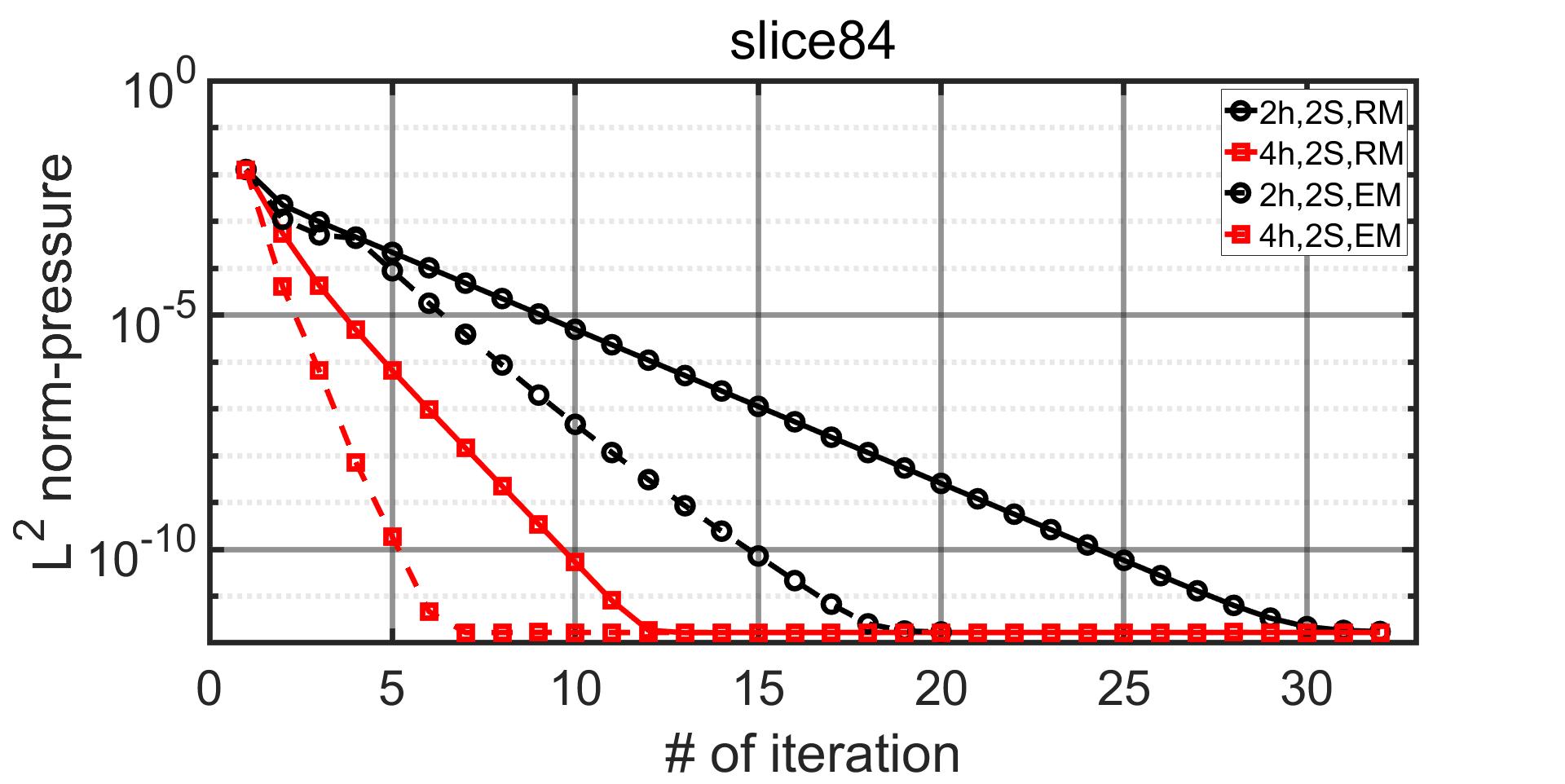}
    \includegraphics[width = 0.48\textwidth]{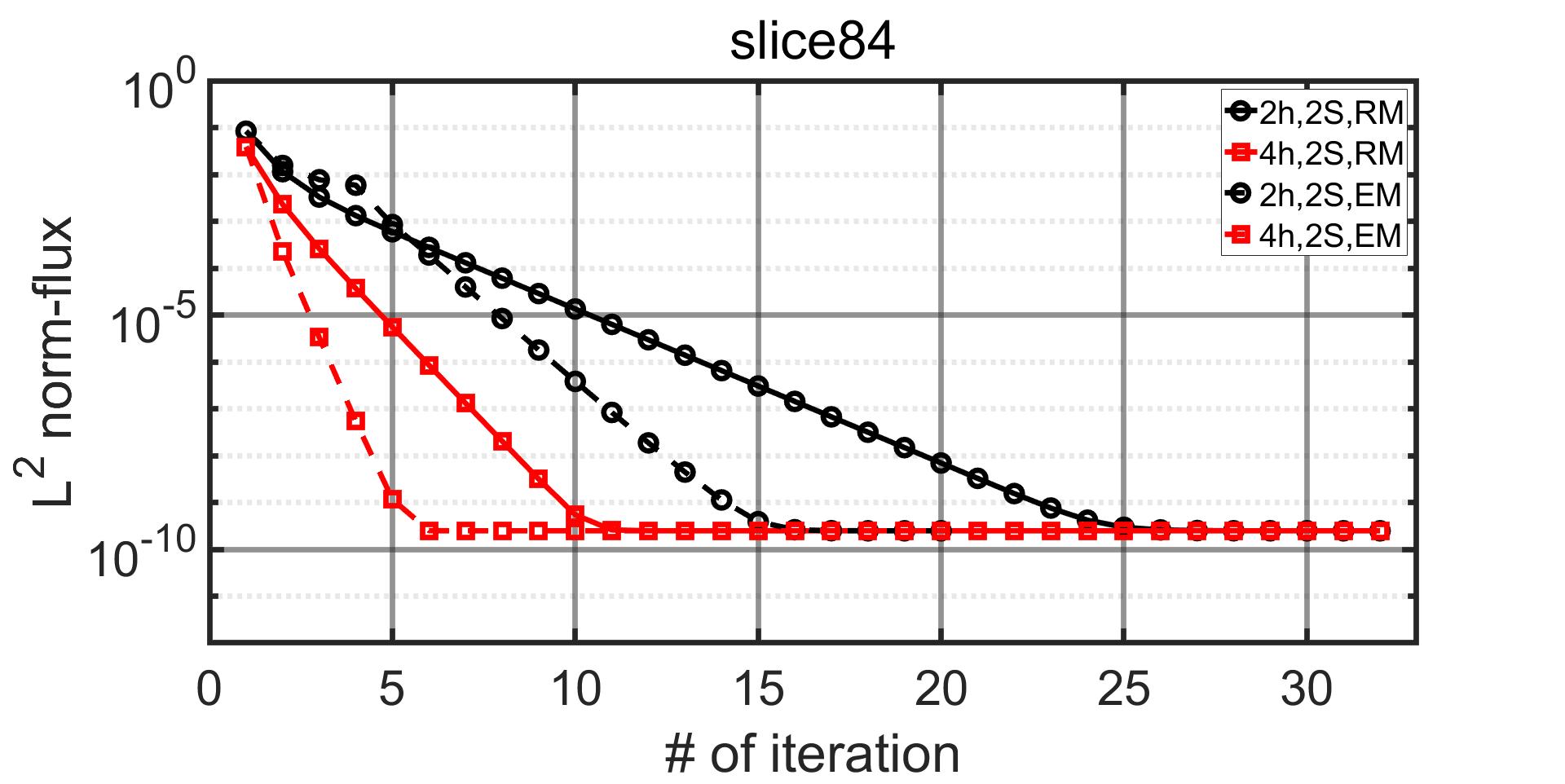}
    \caption{Number of iterations to reach convergence for different oversampling sizes in slice 84: pressure (left) and flux (right).}
    \label{overlap84}
\end{figure}

Here we fix 2 smoothing steps and use oversampling size 2h and 4h. The result is presented in Fig. \ref{overlap42} and \ref{overlap84}.

While increasing oversampling sizes significantly raises the computational cost of computing the set of Multiscale Basis Functions (MBFs), it also greatly reduces the number of iterations required. In slice 42, an additional oversampling of 2h reduces the number of iterations by $30\%$ to $40\%$, while in slice 84, the reduction is approximately $60\%$.

\subsection{The Role Of Smoothing}
Conversely, this section presents the results on how smoothing affects the number of iterations needed for convergence, given fixed oversampling sizes.

\begin{figure}[H]
    \centering
    \includegraphics[width = 0.48\textwidth]{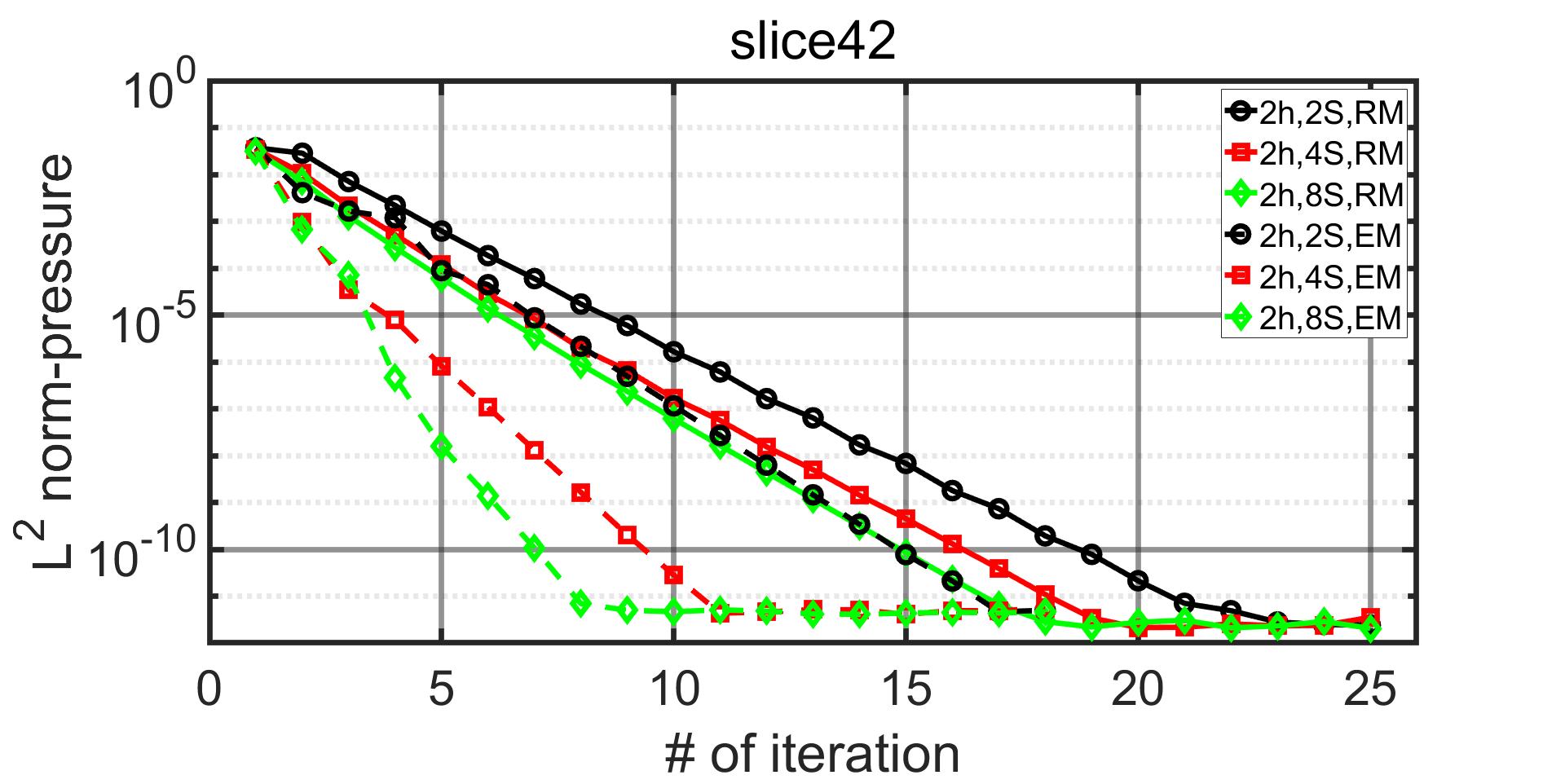}
    \includegraphics[width = 0.48\textwidth]{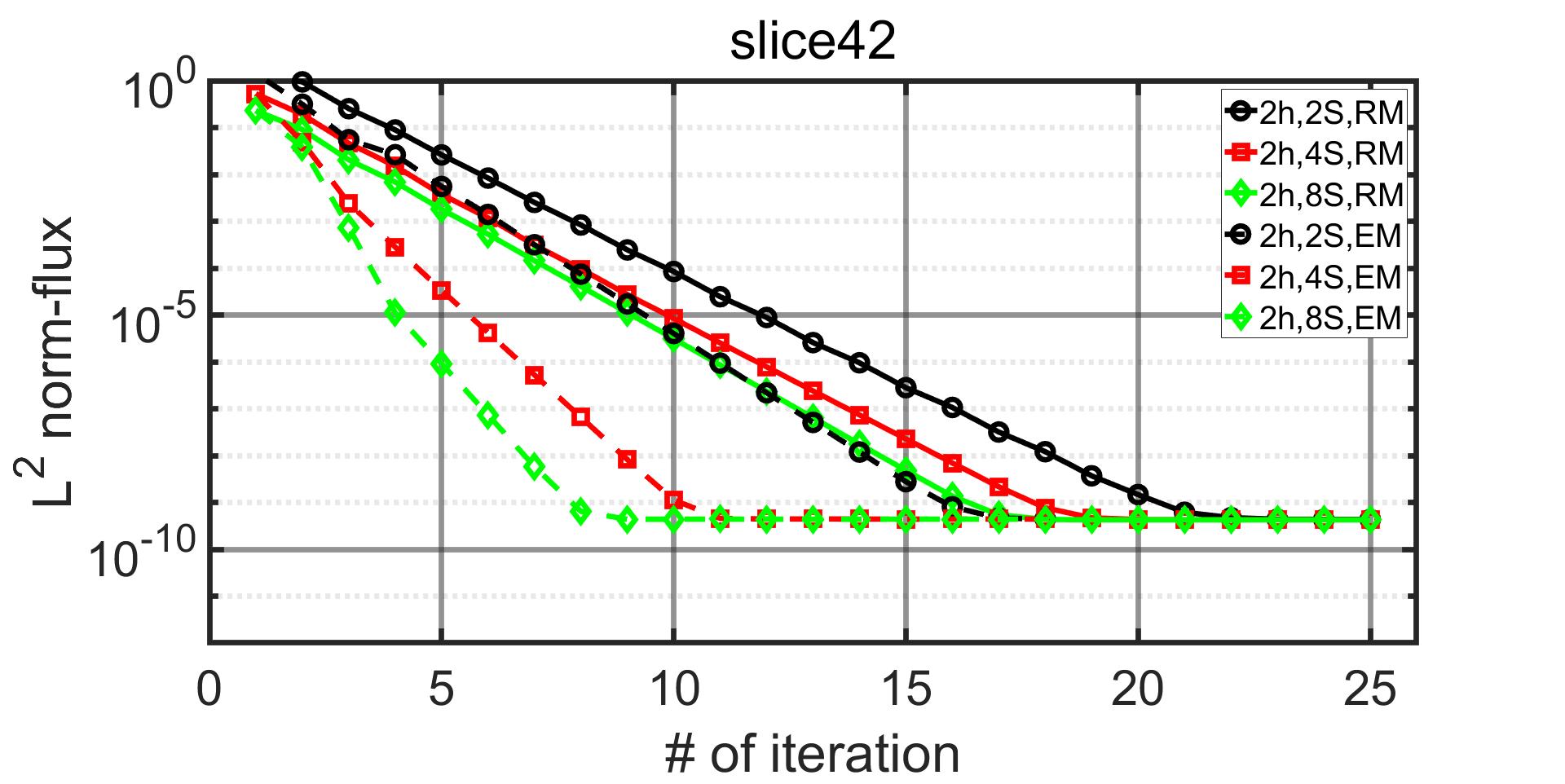}
    \caption{Number of iterations to convergence with varying smoothing steps for slice 42: pressure (left) and flux (right).}
    \label{smoothing42}
\end{figure}

\begin{figure}[H]
    \centering
    \includegraphics[width = 0.48\textwidth]{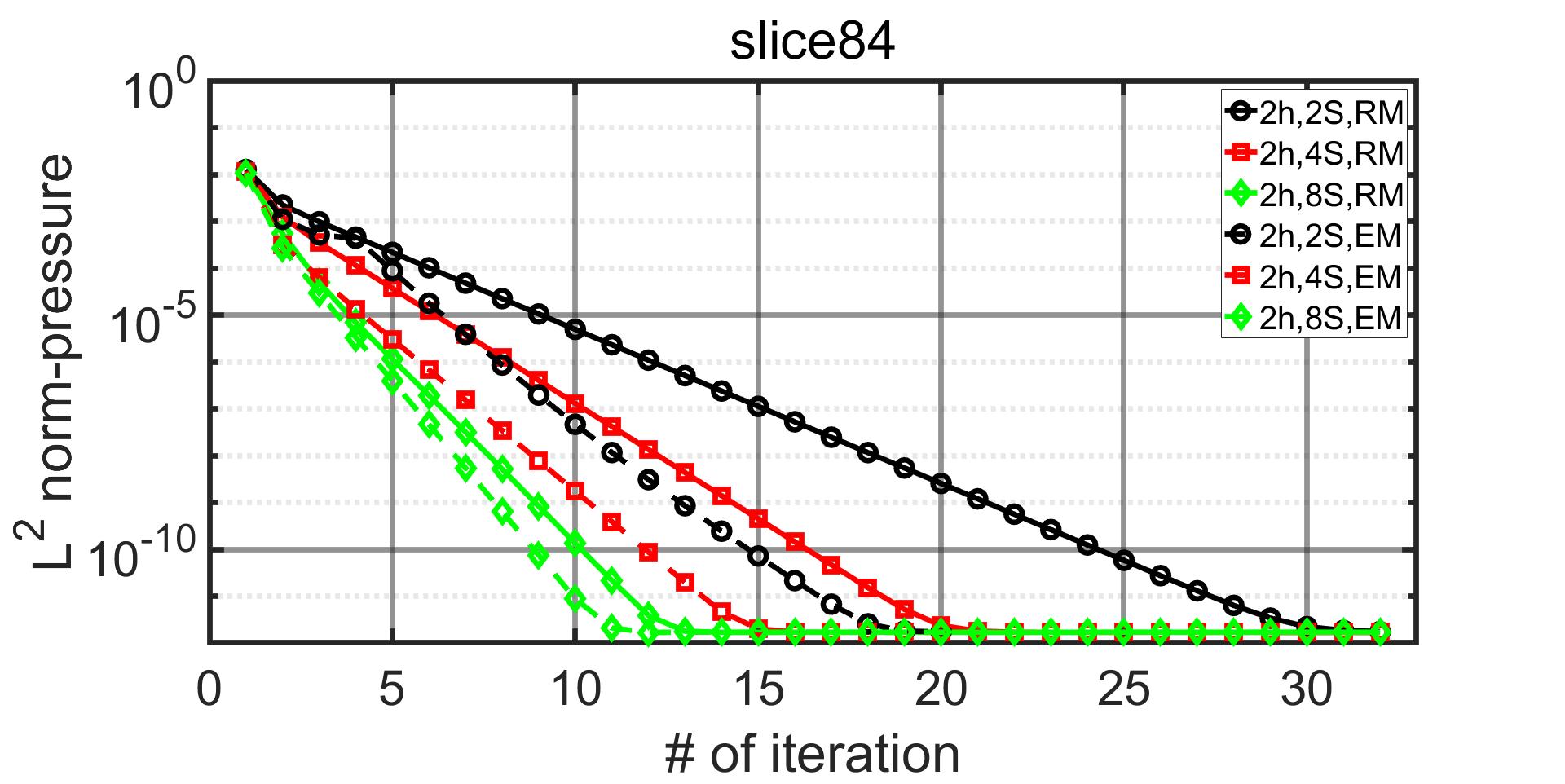}
    \includegraphics[width = 0.48\textwidth]{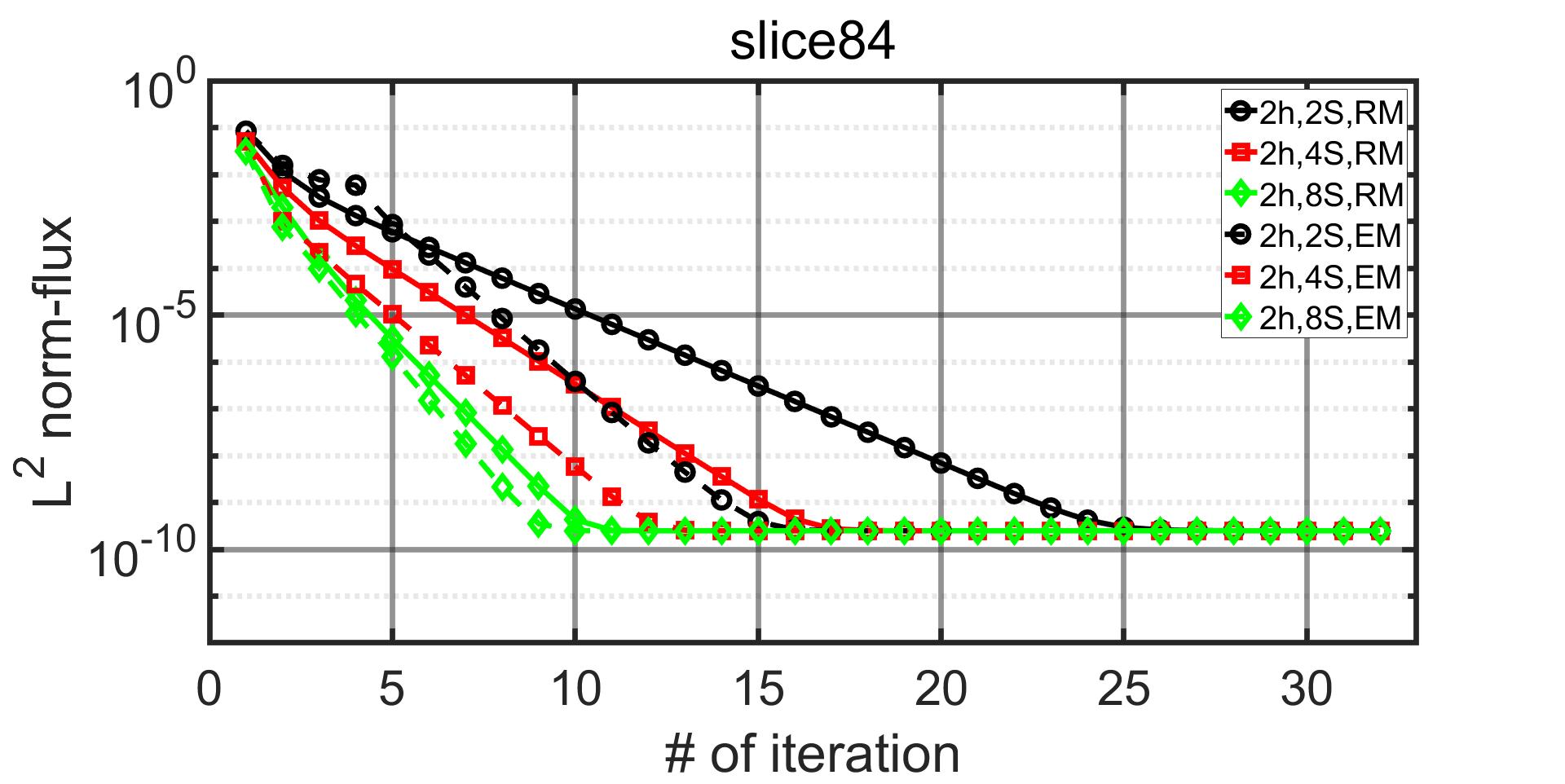}
    \caption{Number of iterations to convergence with varying smoothing steps for slice 84: pressure (left) and flux (right).}
    \label{smoothing84}
\end{figure}

In this analysis, we fix the oversampling size at 2h and use 2, 4, and 8 smoothing steps. The results are shown in Figures \ref{smoothing42} and \ref{smoothing84}.

Increasing the number of smoothing steps also raises computational cost, though not as significantly as in Section 4.3. However, it significantly reduces the number of iterations. Notably, increasing the smoothing steps has a more pronounced impact on slice 42 compared to slice 84, in contrast to the effect of oversampling discussed in Section 4.3.

\subsection{Comparison With An Existing Method}
In that paper \cite{2008FV}, the Iterative Multiscale Finite Volume Method (i-MSFV) is used. The process begins by calculating multiscale basis functions and creating an initial pressure field. During each iteration, smoothing steps are performed with the help of a smoothing operator. The smoothed pressure is then used to compute correction functions, solve the coarse system, and reconstruct the pressure approximation at the end of each iteration.

This paper presents numerous numerical results, from which we select one as an example. The boundary conditions are the same as in our previous section, but the domain decomposition differs. Here, the domain is divided into 4 × 4 subdomains, each containing a local 16 × 16 grid. The source term is set as $q = 1/(h^2)$ in cell (13, 13) and $q = -1/(h^2)$ in cell (32, 32), respectively. The permeability field is assumed to be homogeneous. The paper reports the $L^\infty(\Omega)$ error with respect to the fine grid solution. With 5 smoothing steps per iteration, it reaches an error of $10^{-14}$ after 12 iterations.

\begin{figure}[H]
    \centering
    \includegraphics[width = 0.5\textwidth]{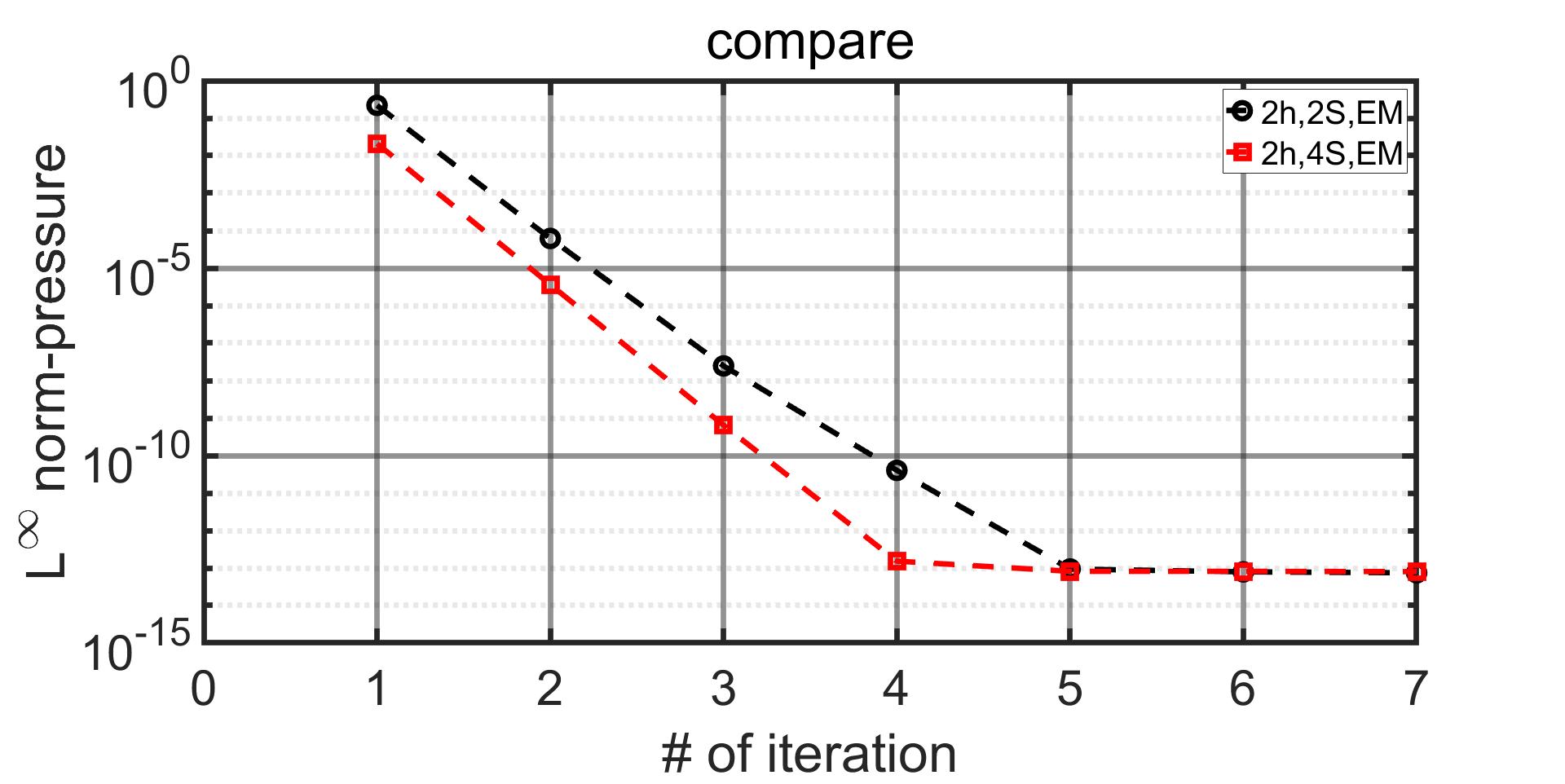}
    \caption{Number of iterations to convergence with varying smoothing steps for the comparison example of the pressure field.}
    \label{compare}
\end{figure}

Figure \ref{compare} presents our results using the extended method for this example. Similar to the result of \ref{compare}, we achieve convergence to an $L^\infty(\Omega)$ error of $10^{-14}$, but more efficiently, in just 4 or 5 iterations using 2 or 4 smoothing steps, respectively. However, comparing computational costs is challenging since we employ different multiscale methods (Multiscale Finite Volume Method and Multiscale Mixed Finite Element Method) and different smoothing procedures.

\section{Conclusions}

In conclusion, this work introduced informed spaces utilizing an offline-online staged method as an extension to MRCM-OS \cite{MRCM-OS}, with the goal of minimizing errors in the pressure field, and subsurface flows while maintaining computational efficiency with very fast converging iterative methods. Two major methods were defined: the Reduced Method (RM) and the Extended Method (EM).

Both RM and EM share the same offline stage. This stage involves calculating the Multiscale Basis Functions (MBFs) using piecewise constant spaces, assembling them into the entire pressure field by solving the interface problems with the same spaces, performing smoothing steps, and finally constructing informed spaces based on the pressure field and subsurface flow along the interfaces of the non-overlapping partitions.

In the online stages, both RM and EM need to update the MBFs according to the informed spaces calculated in the offline stage or from previous iterations. For the interface problem, RM uses the MBFs from the current iteration with piecewise constant spaces, while EM employs MBFs from both the offline stage and the current iteration using linear polynomial spaces. The same smoothing steps are then performed, followed by recalculating the informed spaces as in the offline stage. The iterations continue until the error in the pressure field or subsurface flow meets the specified requirements.

Specifically, our i-MRCM-OS method achieves a flux accuracy of $10^{-10}$ within ten iterations, without requiring large oversampling sizes or numerous smoothing steps. Additionally, comparison with the previous study \cite{2008FV} demonstrates that our method converges more quickly to the same minimum error.

In addition to the two-dimensional numerical examples presented in this paper, informed space also proves beneficial in three-dimensional scenarios. Our forthcoming paper will focus on utilizing informed space to develop a preconditioner, highlighting its potential to create a faster and more accurate preconditioning method.

\bigskip
\noindent {\bf Acknowledgments} 

This material is based upon work supported by the National Science Foundation under Grant No. 2401945. Any opinions, findings, and conclusions or recommendations expressed in this material are those of the authors and do not necessarily reflect the views of the National Science Foundation. The work of F.P. is also partially supported by The University of Texas at Dallas Office of Research and Innovation through the SPARK program.

The authors acknowledge the National Laboratory for Scientific Computing (LNCC/MCTI, Brazil) for providing HPC resources of the S. Dumont supercomputer, which have contributed to the research results reported within this paper (URL: http://sdumont.lncc.br). Additionally, the computing resources of the Cyber-Infrastructure Research Services at the University of Texas at Dallas Office of Information Technology were utilized. The authors would like to express their gratitude to Dr. M\'arcio Borges for his assistance in accessing the LNCC cluster. 

\bigskip
\noindent {\bf Declaration of generative AI and AI-assisted technologies in the writing process}

During the preparation of this work the authors used chatGPT in order to improve language and readability. 
After using this tool, the authors reviewed and edited the content as needed and take full responsibility for the content of the publication.

\bibliography{dilongref1,dilongref2,dilongref3,dilongref4}
\bibliographystyle{unsrt}

\end{document}